\documentclass[12pt]{article}

\usepackage[all]{xy}
\usepackage{amsmath,amssymb,clock, color}
\usepackage[colorlinks, backref]{hyperref}

\newcommand{\Katetov}{Ka\-t\v e\-tov}
\newcommand{\Fraisse}{Fra\"\i\-ss\'e}
\newcommand{\Sierpinski}{Sier\-pi\'n\-ski}

\newcommand{\calC}{\mathcal{C}}
\newcommand{\calF}{\mathcal{F}}
\newcommand{\calP}{\mathcal{P}}
\newcommand{\calR}{\mathcal{R}}

\newcommand{\catA}{\mathcal{A}}
\newcommand{\catC}{\mathcal{C}}

\newcommand{\catV}{\mathcal{V}}

\newcommand{\fin}{{\mathrm{fin}}}
\newcommand{\Ob}{\mathrm{Ob}}
\newcommand{\age}{\mathrm{age}}
\newcommand{\Sym}{\mathrm{Sym}}
\newcommand{\Aut}{\mathrm{Aut}}
\newcommand{\End}{\mathrm{End}}
\newcommand{\id}{\mathrm{id}}
\newcommand{\ID}{\mathrm{ID}}

\newcommand{\NN}{\mathbb{N}}
\newcommand{\QQ}{\mathbb{Q}}
\newcommand{\RR}{\mathbb{R}}

\newcommand{\XYMATRIX}{\xymatrix@M=6pt}
\newcommand{\aremb}{\ar@{^{(}->}}
\newcommand{\dotaremb}{\ar@{^{(}->}^-{\cdot}}
\newcommand{\dothookrightarrow}{\mathrel{\dot\hookrightarrow}}
\newcommand{\pair}[2]{\langle #1, #2\rangle}
\newcommand{\gen}[2]{\langle#1\rangle_{#2}}
\newcommand{\Boxed}[1]{\hbox{$#1$}}

\renewcommand{\sec}{\cap}
\newcommand{\union}{\cup}
\newcommand{\UNION}{\bigcup}
\newcommand{\0}{\emptyset}
\renewcommand{\ge}{\geqslant}
\renewcommand{\le}{\leqslant}
\renewcommand{\phi}{\varphi}
\newcommand{\reducenext}{\underline}

\newcommand{\map}[3]{#1\colon #2 \to #3} 
\newcommand{\cmp}{\circ}

\newcommand{\Z}{\mathbb Z} 
\newcommand{\sn}[1]{\{#1\}}
\newcommand{\koment}[1]{{\color{red}\fbox{\Large#1}}}

\newcommand{\U}{\mathbb U}
\newcommand{\catCh}{\mathcal{C}^{\operatorname{hom}}}
\newcommand{\rest}{\restriction}

\usepackage{theorem}

\theorembodyfont{\slshape}

  \newtheorem{THM}{Theorem}[section]
  \newtheorem{LEM}[THM]{Lemma}
  \newtheorem{PROP}[THM]{Proposition}
  \newtheorem{COR}[THM]{Corollary}

\theorembodyfont{\rmfamily}

  \newtheorem{DEF}[THM]{Definition}
  \newtheorem{EX}[THM]{Example}
  \newtheorem{REM}[THM]{Remark}

\newif\ifQEDsign
\newcommand{\QED}{\global\QEDsigntrue\hfill$\square$}

\newenvironment{PROOF}%
    {\par\noindent\textit{Proof.}\global\QEDsignfalse}%
    {\ifQEDsign\else\QED\fi\par\bigskip\par}

\title{\Katetov\ functors}
\author{Wies\l aw Kubi\'s
\footnote{Research supported by NCN grant 2011/03/B/ST1/00419.}\\
{\small Jan Kochanowski University in Kielce, Poland}\\
{\small Academy of Sciences of the Czech Republic}
\and
Dragan Ma\v sulovi\'c
\footnote{Research supported by the Ministry of Education, Science and Technological Development of the Republic of Serbia, grant 174019.}\\
{\small University of Novi Sad, Serbia}
}

\date{\today \ \clocktime}

\begin{document}
\maketitle

\begin{abstract}
We develop a theory of \emph{\Katetov\ functors} which provide a uniform way of constructing \Fraisse\ limits. Among applications, we present short proofs and improvements of several recent results on the structure of the group of automorphisms and the semigroup of endomorphisms of some \Fraisse\ limits.


\end{abstract}

\tableofcontents

\section{Introduction}\label{kf.sec.introduction}

The theory of \Fraisse\ limits has a long history, inspired by Cantor's theorem saying that the set of rational numbers is the unique, up to isomorphisms, countable linearly ordered set without end-points and such that between any two points there is another one.
In the fifties of the last century Roland \Fraisse\ realized that the ideas behind Cantor's theorem are much more general, and he developed his theory of limits.
Namely, given a class $\catA$ of finitely generated first-order structures with certain natural properties, there exists a unique countably generated structure $L$ (called the \emph{\Fraisse\ limit} of $\catA$) containing isomorphic copies of all structures from $\catA$ and having a very strong homogeneity property, namely, every isomorphism between finitely generated substructures of $L$ extends to an automorphism of $L$.
\Fraisse\ theory can now be called \emph{classical} and it is part of almost every textbook in model theory.
Independently of \Fraisse, around thirty years earlier, Urysohn~\cite{Urysohn} constructed a universal separable complete metric space $\U$ which has the same homogeneity property as \Fraisse\ limits: every finite isometry extends to a bijective isometry of $\U$.
Around the eighties of the last century, merely sixty years after Urysohn's work, \Katetov~\cite{katetov} found a uniform way of extending metric spaces, leading to a new simple construction of the Urysohn space $\U$.
It turns out that \Katetov's construction is functorial, namely, it can be extended to all nonexpansive mappings between metric spaces.

We address the question when a functorial way of constructing a \Fraisse\ limit exists.
Namely, we define the concept of a \emph{\Katetov\ functor} capturing simple extensions of finitely generated structures, whose infinite power gives the \Fraisse\ limit.
The existence of a \Katetov\ functor implies directly that the automorphism group of the \Fraisse\ limit is universal for the class of all automorphism groups of countably generated structures from the given \Fraisse\ class.
  Papers \cite{dolinka-masulovic} and \cite{bilge-melleray} discuss some of the issues addressed in this note,
  without realizing that what one deals with are actually functorial constructions.

As we have mentioned above, our principal motivation comes from \Katetov's construction of the Urys\-ohn space \cite{katetov}, which we briefly recall here in
case of the rational Urysohn space. Let $X$ be a metric space with rational distances. A \emph{\Katetov\ function over $X$}
is every function $\alpha : X \to \QQ$ such that
$$
    |\alpha(x) - \alpha(y)| \le d(x, y) \le \alpha(x) + \alpha(y)
$$
for all $x, y \in X$. Let $K(X)$ be the set of all \Katetov\ functions over $X$.
The $\sup$ metric turns $K(X)$ into a metric space. There is a natural isometric embedding $X \hookrightarrow K(X)$
which takes $a \in X$ to $d(a, \cdot) \in K(X)$. Hence we get a chain of embeddings
$$
  X \hookrightarrow K(X) \hookrightarrow K^2(X) \hookrightarrow K^3(X) \hookrightarrow \cdots
$$
whose colimit is easily seen to be the rational Urysohn space.

It was observed by several authors (see, e.g., \cite{ben-yaacov}, \cite{uspenskij}) that the construction $K$ is actually functorial with respect to embeddings.
Our principal observation is that more is true: if $\catA$ is the category of all finite metric spaces with
rational distances and nonexpansive mappings, and $\catC$ is the category of all countable metric spaces with rational distances and
nonexpansive mappings, then $K$ can be turned into a functor from $\catA$ to $\catC$.
We present the details in the last section.

The paper is organized as follows.
Section~\ref{SecKttvMnDfetc} contains the main concept of a \Katetov\ functor, its basic properties, examples, and a discussion of sufficient conditions for its existence.
We prove, in particular, that a \Katetov\ functor exists if embeddings have pushouts in the category of all homomorphisms.
In Section~\ref{SectKatConstrs} we show how iterations of a \Katetov\ functor lead to \Fraisse\ limits.
It turns out that the \Fraisse\ limit can be viewed as a fixed point of the countable infinite power of a \Katetov\ functor and all orbits of this functor ``tend" to the \Fraisse\ limit, resembling the Banach contraction principle.
Section~\ref{SecBergmanPrpty} deals with the semigroup Bergman property.
We prove that in the presence of a \Katetov\ functor, under some mild additional assumptions the endomorphism monoid $\End(L)$ of the \Fraisse\ limit $L$ is strongly distorted and its Sierpi\'nski rank is at most five.
Applying a result from \cite{maltcev-mitchell-ruskuc}, we conclude that if $\End(L)$ is not finitely generated, then it has the Bergman property.
This extends a recent result of Dolinka~\cite{dolinka}.
The last Section~\ref{SecAppndxKattv} is an appendix containing description of the original \Katetov\ functor on metric spaces with nonexpansive mappings.

\subsection{The setup}

Let $\Delta = \calR \union \calF \union \calC$ be a first-order language, where $\calR$ is a set of relational symbols,
$\calF$ a set of functional symbols, and $\calC$ a set of constant symbols. We say that $\Delta$ is
a \emph{purely relational language} if $\calF = \calC = \0$. For a $\Delta$-structure $A$ and $X \subseteq A$, by
$\gen XA$ we denote the substructure of $A$ generated by $X$. We say that $A$ is \emph{finitely generated} if
$A = \gen XA$ for some finite $X \subseteq A$.
The fact that $A$ is a substructure of $B$ will be denoted by $A \le B$.

Let $\catC$ be a category of $\Delta$-structures.
A \emph{chain in $\catC$} is a sequence of objects and embeddings of the form
$
  C_1 \hookrightarrow C_2 \hookrightarrow C_3 \hookrightarrow \cdots
$.
Note that although there may be other kinds of morphisms in $\catC$, a chain always consists of objects and embeddings.
We shall say that $L$ is a \emph{standard colimit} of the chain $C_1 \hookrightarrow C_2 \hookrightarrow \cdots$ if it is a colimit of this chain in the usual sense and moreover, after forgetting the structure $L$ is still a colimit in the category of sets.
In other words, if the embeddings are inclusions, that is, $C_1 \le C_2 \le \cdots$ then a standard colimit is $L = \bigcup_{n\in \NN}C_n$ with an appropriate $\Delta$-structure making it a colimit in $\catC$.
We shall say that $\catC$ has \emph{standard colimits of chains} if every chain in $\catC$ has a standard colimit in $\catC$.
Given $C \in \Ob(\catC)$, let $\Aut(C)$ denote the permutation group consisting of all automorphisms of $C$, and let $\End(C)$
denote the transformation monoid consisting of all $\catC$-morphisms $C \to C$.
It may be the case that $\End(C)$ consists of all embeddings of $C$ into $C$ (if $\catC$ consists of embeddings only).
We shall sometimes write $\End_\catC(C)$ instead of $\End(C)$ in order to emphasize that we consider $\catC$-morphisms only.
Let $\age(C)$ denote the class of all
finitely generated objects that embed into $C$.
We say that $\catA$ has the \emph{joint-embedding property} (briefly: (JEP)) if every two structures in $\catA$ embed into a common structure in $\catA$.

\paragraph{Standing assumption.}
Throughout the paper we assume the following. Let $\Delta$ be a first-order language,
let $\catC$ be a category of countably generated $\Delta$-structures
and some appropriately chosen class of morphisms that includes all embeddings (and hence all isomorphisms).
Let $\catA$ be the full subcategory of $\catC$ spanned by all finitely generated structures in $\catC$.
In particular, $\catA$ is \emph{hereditary} in the sense that given $A \in \Ob(\catA)$, every finitely generated substructure\footnote{Recall that substructures of finitely generated structures may not be finitely generated. For example, the free group with 2 generators has a subgroup isomorphic to the free group with infinitely many generators.} of $A$ is an object of $\catA$.

We assume that the following holds:
\begin{itemize}
\item
  $\catC$ has standard colimits of chains;
\item
  every $C \in \Ob(\catC)$ is a colimit of some chain $A_1 \hookrightarrow A_2 \hookrightarrow \cdots$ in $\catA$;
\item
  $\catA$ has only countably many isomorphism types; and
\item
  $\catA$ has the joint embedding property (JEP).
\end{itemize}

\bigskip

We say that $C \in \Ob(\catC)$ is a \emph{one-point extension} of $B \in \Ob(\catC)$ if there is an embedding
$j : B \hookrightarrow C$ and an $x \in C \setminus j(B)$ such that $C = \gen{j(B) \union \{x\}}{C}$.
In that case we write $j : B \dothookrightarrow C$ or simply $B \dothookrightarrow C$.

The following lemmas are immediate consequences of the fact that $\catC$ is a category of $\Delta$-structures
and the fact that $\catA$ is spanned by finitely generated objects in $\catC$.

\begin{LEM}[Reachability]\label{kf.lem.reach}
  $(a)$ For all $A, B \in \Ob(\catA)$ and an embedding $A \hookrightarrow B$ which is not an isomorphism,
  there exist an $n \in \NN$ and $A_1, \ldots, A_n \in \Ob(\catA)$ such that
  $
    A \dothookrightarrow  A_1 \dothookrightarrow A_2 \dothookrightarrow \cdots \dothookrightarrow A_n = B.
  $

  $(b)$ For all $C, D \in \Ob(\catC)$ and an embedding $f : C \hookrightarrow D$ which is not an isomorphism,
  there exist $C_1, C_2 \ldots \in \Ob(\catC)$ such that
  $$
    \XYMATRIX{
      C \dotaremb[r] \aremb[drr]_-f & C_1 \dotaremb[r] \aremb[dr] & C_2 \dotaremb[r] \aremb[d] & \cdots \aremb[dl]\\
                     &                  &                  D
    }
  $$
  is a colimit diagram in $\catC$.
\end{LEM}

\begin{LEM}\label{kf.lem.c1pe}
  Let $C, D \in \Ob(\catC)$ be structures such that $f : C \dothookrightarrow D$ and let
  $A_1 \hookrightarrow A_2 \hookrightarrow \ldots$ be a chain in $\catA$ whose colimit is~$C$.
  Then there exists a chain $B_1 \hookrightarrow B_2 \hookrightarrow \ldots$ in $\catA$ whose colimit is~$D$ and the following diagram commutes
  $$
    \XYMATRIX{
      & \\
      A_1 \aremb[r] \dotaremb[d] \aremb@/^9mm/[rrrr] & A_2 \aremb[r] \dotaremb[d] \aremb@/^5mm/[rrr] & A_3 \aremb[r] \dotaremb[d] \aremb@/^2mm/[rr] & \cdots & C \dotaremb[d]_-f\\
      B_1 \aremb[r] \aremb@/_9mm/[rrrr] & B_2 \aremb[r] \aremb@/_5mm/[rrr] & B_3 \aremb[r] \aremb@/_2mm/[rr]& \cdots & D\\
      & \\
    }
  $$
where the curvy arrows are canonical embeddings into the colimits.
\end{LEM}
\begin{PROOF}
  Without loss of generality we can assume that $C \le D$, and that $A_1 \le A_2 \le \ldots \le C$, so that
  $C = \UNION_{i \in \NN} A_i$.
  Since $D$ is a one-point extension of $C$, there exists an $x \in D \setminus C$ such that $D = \langle C \union \{x\} \rangle_D$.
  Put $B_i = \langle A_i \union \{x\} \rangle_D$.
\end{PROOF}

The next lemma is rather obvious, as we assume that colimits are standard.

\begin{LEM}[Factoring through the colimit of a chain]\label{kf.lem.factoring-colim}
  Let

  $$C_1 \hookrightarrow C_2 \hookrightarrow \cdots$$
  be a chain in $\catC$ and let $L$ be its colimit
  with the canonical embeddings $\iota_k : C_k \hookrightarrow L$.
  Then for every $A \in \Ob(\catA)$ and every morphism $f : A \to L$ there is an $n \in \NN$ and a
  morphism $g : A \to C_n$ such that $f \circ g = \iota_n$.
  Moreover, if $f$ is an embedding, then so is $g$.
  $$
    \XYMATRIX{
        & C_n \aremb[d]^-{\iota_n}\\
        A \ar[ur]^-g \ar[r]_-f & L
    }
  $$
\end{LEM}

\begin{LEM}\label{kf.lem.age}
  For every $C \in \Ob(\catC)$ we have that $\age(C) \subseteq \Ob(\catA)$.
\end{LEM}
\begin{PROOF}
  Take any $C \in \Ob(\catC)$, and let $A_1 \hookrightarrow A_2 \hookrightarrow \cdots$ be a chain in $\catA$
  whose colimit is $C$. Take any $B \in \age(C)$. Then $B \hookrightarrow C$, so by
  Lemma~\ref{kf.lem.factoring-colim} there is an $n \in \NN$ and an
  embedding $g : B \hookrightarrow A_n$ such that
  $$
    \XYMATRIX{
        & A_n \aremb[d]\\
        B \aremb[ur]^-g \aremb[r] & C
    }
  $$
  Therefore, $B \hookrightarrow A_n \in \Ob(\catA)$, so the assumption that $\catA$ is hereditary yields
  $B \in \Ob(\catA)$.
\end{PROOF}

\section{\Katetov\ functors}\label{SecKttvMnDfetc}

\begin{DEF}\label{kf.def.KatetovFtr}
  A functor $K^0 : \catA \to \catC$ is a \emph{\Katetov\ functor} if:
  \begin{itemize}
  \item
    $K^0$ preserves embeddings, that is, if $f : A \rightarrow B$ is an embedding in $\catA$, then
    $K^0(f) : K^0(A) \rightarrow K^0(B)$ is an embedding in $\catC$; and
  \item
    there is a natural transformation $\eta^0 : \ID \to K^0$ such that
    for every one-point extension $A \dothookrightarrow B$ where $A, B \in \Ob(\catA)$, there is an embedding
    $g : B \hookrightarrow K^0(A)$ satisfying
    \begin{equation}\label{kf.eq.def}
      \XYMATRIX{
        A \aremb[r]^-{\eta^0_A} \dotaremb[d] & K^0(A) \\
        B \aremb[ur]_-{g}
      }
    \end{equation}
  \end{itemize}
\end{DEF}

\begin{THM}\label{kf.thm.ext-to-C}
  If there exists a \Katetov\ functor $K^0 : \catA \to \catC$ then there is a functor $K : \catC \to \catC$ such that:
  \begin{itemize}
  \item
    $K$ is an extension of $K^0$ (that is, $K$ and $K^0$ coincide on $\catA$);
  \item
    there is a natural transformation $\eta : \ID \to K$ which is an extension of $\eta^0$
    (that is, $\eta_A = \eta^0_A$ whenever $A \in \Ob(\catA)$);
  \item
    $K$ preserves embeddings.
  \end{itemize}
\end{THM}
\begin{PROOF}
  The obvious candidate for $K$ is the left Kan extension of $K^0$ along the inclusion functor $E : \catA \to \catC$
  (which acts identically on both objects and morphisms of $\catA$).
  To show that such an extension exists it suffices to show that the diagram
  $(E \downarrow C) \overset{\Pi}\longrightarrow \catA \overset{K^0}\longrightarrow \catC$
  has a colimit in $\catC$ for every $C \in \Ob(\catC)$, where $\Pi$ is the projection functor from the
  comma category $(E \downarrow C)$ to $\catA$ which takes an object $(A, h : A \to C)$ of the comma category to
  its first coordinate $A$, and acts on morphisms accordingly~\cite{MacLane}.

  Take any $C \in \Ob(\catC)$ and let $A^C_1 \hookrightarrow A^C_2 \hookrightarrow \cdots$ be a chain in $\catA$ whose
  colimit is $C$. Let $\iota^C_n : A^C_n \hookrightarrow C$ be the canonical embeddings.
  Recall that for every $B \in \Ob(\catA)$ and every morphism $f : B \to C$ there is an $n$ and a morphism
  $f_n : B \to A^C_n$ such that $\iota^C_n \circ f_n = f$ (Lemma~\ref{kf.lem.factoring-colim}):
  $$
    \XYMATRIX{
        A^C_n \aremb[dr]^-{\iota^C_n} \aremb[r] & A^C_{n+1} \aremb[d]^-{\iota^C_{n+1}}\\
        B \ar[u]^-{f_n} \ar[r]_-f & C
    }
  $$
  The diagram $(E \downarrow C) \overset{\Pi}\longrightarrow \catA \overset{K^0}\longrightarrow \catC$
  then takes the form
  $$
    \XYMATRIX{
        K^0(A^C_n) \aremb[r] & K^0(A^C_{n+1}) \\
        K^0(B) \ar[u]^-{K^0(f_n)}
    }
  $$
  Let $D$ be the colimit of the chain $K^0(A^C_1) \hookrightarrow K^0(A^C_2) \hookrightarrow \cdots$ with the
  canonical embeddings $\iota^D_n : K^0(A^C_n) \hookrightarrow D$. For each $f : B \to C$ in $\catA$ let
  $f' = \iota^D_n \circ K_0(f_n) : K^0(B) \to D$:
  $$
    \XYMATRIX{
        K^0(A^C_n) \aremb[dr]^-{\iota^D_n} \aremb[r] & K^0(A^C_{n+1}) \aremb[d]^-{\iota^D_{n+1}}\\
        K^0(B) \ar[u]^-{K^0(f_n)} \ar[r]_-{f'} & D
    }
  $$
  Then it is easy to show that $D$ is the colimit of the diagram
  $(E \downarrow C) \overset{\Pi}\longrightarrow \catA \overset{K^0}\longrightarrow \catC$ in $\catC$.
  Therefore, $K^0$ has the left Kan extension $K$ along $E$.

  Let us show that $K$ preserves embeddings. Take any embedding $f : C \hookrightarrow D$ in $\catC$.
  Let $A^C_1 \hookrightarrow A^C_2 \hookrightarrow \cdots$
  be a chain in $\catA$ whose colimit is $C$ and let $A^D_1 \hookrightarrow A^D_2 \hookrightarrow \cdots$
  be a chain in $\catA$ whose colimit is $D$. Moreover,
  let $\iota^C_n : A^C_n \hookrightarrow C$ and $\iota^D_n : A^D_n \hookrightarrow D$ be the corresponding
  canonical embeddings. By Lemma~\ref{kf.lem.factoring-colim}, for every $k$ there is an $n_k$ and
  a morphism $f_k : A^C_k \to A^D_{n_k}$ (which is necessarily an embedding) such that
  $$
    \XYMATRIX{
        A^C_k \aremb[r]^-{\iota^C_k} \aremb[d]_-{f_k} & C \aremb[d]^-f\\
        A^D_{n_k} \aremb[r]_-{\iota^D_{n_k}} & D
    }
  $$
  Without loss of generality $n_k$'s can be chosen in such a way that $n_1 < n_2 < \ldots$. In the extension we then have
  $$
    \XYMATRIX{
        K^0(A^C_k) \aremb[r]^-{\iota^{K(C)}_k} \aremb[d]_-{K^0(f_k)} & K(C) \ar[d]^-{K(f)}\\
        K^0(A^D_{n_k}) \aremb[r]_-{\iota^{K(D)}_{n_k}} & K(D)
    }
  $$
  whence follows that $K(f)$ is also an embedding.

  Analogous argument provides a construction of the natural transformation $\eta : \ID \to K$
  which extends $\eta^0$. Consider the diagram:
  $$
    \XYMATRIX{
        A^C_{n} \aremb@/^5mm/[rr]^-{\iota^C_{n}} \aremb[d]_-{\eta^0_{A^C_{n}}} \aremb[r]
        &  A^C_{n+1} \aremb[r]_-{\iota^C_{n+1}} \aremb[d]_-{\eta^0_{A^C_{n+1}}}
        & C \ar[d]^-{\eta_C}
        \\
        K^0(A^C_{n}) \aremb@/_5mm/[rr]_-{\iota^{K(C)}_{n}} \aremb[r]
        & K^0(A^C_{n+1}) \aremb[r]^-{\iota^{K(C)}_{n+1}}
        & K(C)
    }
  $$
  Since $C$ is the colimit of the chain $A^C_1 \hookrightarrow A^C_2 \hookrightarrow \cdots$, there is a unique
  morphism $\eta_C : C \to K(C)$ to the tip of the competing compatible cone. The morphism $\eta_C$ is
  clearly an embedding and it is easy to check that all the morphisms $\eta_C$ constitute a natural
  transformation $\eta : \ID \to K$.
\end{PROOF}

We also say that $K$ is a \emph{\Katetov\ functor} and from now on we denote both $K$ and $K^0$ by~$K$, and both $\eta$ and $\eta^0$ by~$\eta$.
An obvious yet important property of $K$ is that all its powers remain \Katetov.
Specifically, for $n \in \NN$ define $\eta^n : \ID \to K^n$ as
$\eta^n_C = \eta_{K^{n-1}(C)} \circ \ldots \circ \eta_{K(C)} \circ \eta_C : C \to K^n(C)$.
Then $\eta^n$ is a natural transformation witnessing that $K^n$ is a \Katetov\ functor.
We shall elaborate this in Section~\ref{SectKatConstrs}.
For now, we state the following important property of finite iterations of $K$.

\begin{LEM}\label{kf.lem.n-chain}
  Let $K : \catA \to \catC$ be a \Katetov\ functor. Then for every embedding $g : A \hookrightarrow B$,
  where $A, B \in \Ob(\catA)$, there is an $n \in \NN$ and an embedding
  $h : B \hookrightarrow K^n(A)$ satisfying $h \circ g = \eta_A^n$.
  $$
      \XYMATRIX{
        A \aremb[r]^-{\eta^n_A} \aremb[d]_{g} & K^n(A) \\
        B \aremb[ur]_{h}
      }
  $$
\end{LEM}
\begin{PROOF}
  If $g$ is an isomorphism, take $n = 1$ and $h = \eta_A \circ g^{-1}$. Assume, therefore, that
  $g$ is not an isomorphism. Then by Lemma~\ref{kf.lem.reach}$(a)$ there exist $n \in \NN$
  and $A_1, \ldots, A_n \in \Ob(\catA)$ such that
  $$
    A \dothookrightarrow A_1 \dothookrightarrow A_2 \dothookrightarrow \cdots \dothookrightarrow A_n = B.
  $$
  It is easy to see that the diagram in Fig.~\ref{kf.fig.1} commutes: the triangles commute by the definition
  of a \Katetov\ functor, while the parallelograms commute because $\eta$ is a natural transformation.
  So, take $h = K^{n-1}(f_1) \circ K^{n-2}(f_2) \circ \ldots \circ K(f_{n-1}) \circ f_n$.
\end{PROOF}

\begin{figure}
  $$
    \XYMATRIX{
      A   \dotaremb[d] \aremb[dr]^{\eta_A}\\
      A_1 \dotaremb[d] \aremb[dr]^{\eta_{A_1}} \aremb[r]^-{f_1} & K(A) \aremb[dr]^{\eta_{K(A)}} \\
      A_2 \dotaremb[d] \aremb[dr]^{\eta_{A_2}} \aremb[r]^-{f_2} & K(A_1) \aremb[dr]^{\eta_{K(A_1)}} \aremb[r]^{K(f_1)} & K^2(A) \aremb[dr]^{\eta_{K^2(A)}}\\
      A_3 \dotaremb[d] \aremb[dr]^{\eta_{A_3}} \aremb[r]^-{f_3} & K(A_2) \aremb[dr]^{\eta_{K(A_2)}} \aremb[r]^{K(f_2)} & K^2(A_1) \aremb[dr]^{\eta_{K^2(A_1)}} \aremb[r]^{K^2(f_1)} & K^3(A) \aremb[dr]^{\eta_{K^3(A)}} \\
      {\vdots} \dotaremb[d] \aremb[dr]^{\eta_{A_{n-1}}} & {\ddots} \aremb[dr]^{\eta_{K(A_{n-2})}} & {\ddots} \aremb[dr]^{\eta_{K^2(A_{n-3})}} & {\ddots} \aremb[dr]^{\eta_{K^{n-2}(A_1)}} & {\ddots} \aremb[dr]^{\eta_{K^{n-1}(A)}}\\
      A_n \aremb[r]_-{f_n} & K(A_{n-1}) \aremb[r]_{K(f_{n-1})} & K^2(A_{n-2}) \aremb[r]_-{K^2(f_{n-2})} & \cdots \aremb[r]_-{K^{n-2}(f_2)} & K^{n-1}(A_1) \aremb[r]_-{K^{n-1}(f_1)} & K^n(A)
    }
  $$
\caption{The proof of Lemma~\ref{kf.lem.n-chain}}
\label{kf.fig.1}
\end{figure}

\subsection{Examples}\label{kf.subsec.ex}

Below we collect some examples of \Katetov\ functors.
The second one shows in particular that \Katetov\ functors (as well as their powers) do not necessarily have the extension property for embeddings between objects of $\catC$.
In other words, Lemma~\ref{kf.lem.n-chain} does not hold for embeddings between $\catC$-objects.

\begin{EX}
  \textit{A \Katetov\ functor on the category of finite metric spaces with rational distances and nonexpansive mappings}

This is a small modification of the original \Katetov\ functor, in order to fit into \Fraisse\ theory.
The details are explained in Section~\ref{SecAppndxKattv} below.
\end{EX}

Let $\calP_2(X) = \{Y \subseteq X : |Y| = 2\}$, and let $\calP_\fin(X)$ denote the set of all finite subsets of~$X$.

\begin{EX}
  \textit{A \Katetov\ functor on the category of graphs and graph homomorphisms.}
  Let $\pair V E$ be a graph, where $E \subseteq \calP_2(V)$. Put $K(\pair V E) = \pair{V^*}{E^*}$ where
  \begin{align*}
    V^* &= V \union \calP_\fin(V), \\
    E^* &= E \union \{ \{v, A\} : A \in \calP_\fin(V), v \in A \}.
  \end{align*}
  For a graph homomorphism $f : \pair{V_1}{E_1} \to \pair{V_2}{E_2}$ let $f^* = K(f)$ be a mapping from
  $V_1^*$ to $V_2^*$ defined by $f^*(v) = f(v)$ for $v \in V_1$ and $f^*(A) = f(A)$ for $A \in \calP_\fin(V_1)$.
  Then it is easy to show that $f^*$ is a graph homomorphism from $\pair{V_1^*}{E_1^*}$ to $\pair{V_2^*}{E_2^*}$.
  Moreover, if $f$ is an embedding, then so is $f^*$.

Now let $G$ be an infinite graph.
Let $H = G \cup \{v\}$, where $v$ is connected to all the vertices of $G$.
Note that each vertex of $K(G) \setminus G$ has a finite degree in $K(G)$.
Thus, there is no embedding of $H$ extending $\eta_G \colon G \to K(G)$.
This shows that $K$ does not have the extension property for embeddings in the bigger category consisting of all countable graphs.
The same holds for $K^n$ for every $n>1$ (and even for its $\omega$-power), because all ``new'' vertices in $K^n(G)$ have finite degrees in $G$.
\end{EX}

\begin{EX}\label{kf.ex.Kn-free}
  \textit{A \Katetov\ functor on the category of $K_n$-free graphs and graph embeddings.}
  Fix an integer $n \ge 3$. Let $\pair V E$ be a $K_n$-free graph, where $E$ is the set of some 2-element subsets of~$V$.
  Put $K(\pair VE) = \pair{V^*}{E^*}$ where
  \begin{align*}
    V^* = V & \mathstrut \union V', \\
            & V' = \{A \in \calP_\fin(V) : \pair{A}{E \sec \calP_2(A)} \text{ is $K_{n-1}$-free}\}, \\
    E^* = E & \mathstrut \union \{ \{v, A\} : A \in V', v \in A \}.
  \end{align*}
  For a graph embedding $f : \pair{V_1}{E_1} \hookrightarrow \pair{V_2}{E_2}$ let $f^* = K(f)$ be a mapping from
  $V_1^*$ to $V_2^*$ defined by $f^*(v) = f(v)$ for $v \in V_1$ and $f^*(A) = f(A)$ for $A \in V_1'$.
  Then it is easy to show that $f^*$ is a graph embedding from $\pair{V_1^*}{E_1^*}$ to $\pair{V_2^*}{E_2^*}$.
\end{EX}

\begin{EX}
  \textit{A \Katetov\ functor on the category of digraphs and digraph homomorphisms.}
  Let $\pair V E$ be a digraph, where $E \subseteq V^2$ is an irreflexive relation satisfying $(x, y) \in E \Rightarrow (y, x) \notin E$.
  Put $K(\pair V E) = \pair{V^*}{E^*}$ where
  \begin{align*}
    V^* = V & \mathstrut \union V', \\
            & V' = \{\pair AB :  A, B \in \calP_\fin(V) \text{ such that } A \sec B = \0\}, \\
    E^* = E & \mathstrut \union \{ \pair{v}{\pair AB} : v \in V, \pair AB \in V', v \in A \}\\
            & \mathstrut \union \{ \pair{\pair AB}{v} : v \in V, \pair AB \in V', v \in B \}.
  \end{align*}
  For a digraph homomorphism $f : \pair{V_1}{E_1} \to \pair{V_2}{E_2}$ let $f^* = K(f)$ be a mapping from
  $V_1^*$ to $V_2^*$ defined by $f^*(v) = f(v)$ for $v \in V_1$ and $f^*(\pair AB) = \pair{f(A)}{f(B)}$ for $\pair AB \in V_1'$.
  Then it is easy to show that $f^*$ is a digraph homomorphism from $\pair{V_1^*}{E_1^*}$ to $\pair{V_2^*}{E_2^*}$.
  Moreover, if $f$ is an embedding, then so is $f^*$.
\end{EX}

\begin{EX}
  \textit{A \Katetov\ functor on the category of all finite linear orders and monotone mappings.}
  For a linear order $\pair{A}{\Boxed{\le}}$ put $K(\pair{A}{\Boxed{\le}}) = \pair{A^*}{\Boxed{\le^*}}$ where
  \begin{align*}
    A^* = A &\mathstrut\union A', \\
            &A' = \{\pair UV : \{U, V\} \text{ is a partition of $A$ and } \forall u \in U \; \forall v \in V\;(u \le v)\},\\
    \Boxed{\le^*} = \Boxed{\le} &\mathstrut\union \{ \pair{a}{\pair UV} : a \in U\} \union \{ \pair{\pair UV}{a} : a \in V\}\\
                                &\mathstrut\union \{ \pair{\pair{U_1}{V_1}}{\pair{U_2}{V_2}} : V_1 \sec U_2 \ne \0\}\\
                                &\mathstrut\union \{ \pair{\pair UV}{\pair UV} : \pair UV \in A'\}.
  \end{align*}
  Then it is easy to see that $\le^*$ is a linear order on $A^*$.
  For a monotone map $f \colon \pair{A_1}{\Boxed{\le_1}} \to \pair{A_2}{\Boxed{\le_2}}$ let $f^* = K(f)$ be the mapping from
  $A_1^*$ to $A_2^*$ defined by $f^*(a) = f(a)$ for $a \in A_1$ and for $\pair UV \in A_1'$
  we put $f^*(\pair UV) = \pair{A_2 \setminus W}{W}$ where $W = f(V)$.
  Then it is easy to show that $f^*$ is a monotone map from $\pair{A_1^*}{\Boxed{\le_1^*}}$ to $\pair{A_2^*}{\Boxed{\le_2^*}}$.
  Moreover, if $f$ is an embedding, then so is $f^*$.

  The description of $K(\pair{A}{\Boxed{\le}})$ in case $\pair{A}{\Boxed{\le}}$ is a countable linear order is more involved.
  As an illustration, let us just say that $K(\pair{\QQ}{\Boxed{\le}})$ is of the form $Q_1 \cup Q_2$,
  where both $Q_1$ and $Q_2$ are dense in $K(\pair{\QQ}{\Boxed{\le}})$ and $Q_1$ serves as a copy of $\pair{\QQ}{\Boxed{\le}}$ while
  $Q_2$ serves as the set of all one-point extensions of finite subsets of $Q_1$, ordered in a suitable way.

\end{EX}

\begin{EX}
  \textit{A \Katetov\ functor on the category of partially ordered sets and monotone mappings.}
  For a partially ordered set $\pair{A}{\Boxed{\le}}$
  put $K(\pair{A}{\Boxed{\le}}) = \pair{A^*}{\Boxed{\le^*}}$ where
  \begin{align*}
    A^* = A &\mathstrut\union A',\\
            &A' = \{\pair UV : U, V \in \calP_\fin(A) \text{ and } \forall u \in U \; \forall v \in V\;(u \le v)\}, \\
    \Boxed{\le^*} = \Boxed{\le} &\mathstrut\union \{ \pair{a}{\pair UV} : \exists u \in U \; (a \le u)\}\\
                                &\mathstrut\union \{ \pair{\pair UV}{a} : \exists v \in V \; (v \le a)\}\\
                                &\mathstrut\union \{ \pair{\pair{U_1}{V_1}}{\pair{U_2}{V_2}} : \exists v \in V_1 \; \exists u \in U_2 \; (v \le u)\}\\
                                &\mathstrut\union \{ \pair{\pair UV}{\pair UV} : \pair UV \in A'\}.
  \end{align*}
  Then it is easy to see that $\le^*$ is a partial order on $A^*$.
  For a monotone map $f : \pair{A_1}{\Boxed{\le_1}} \to \pair{A_2}{\Boxed{\le_2}}$ let $f^* = K(f)$ be a mapping from
  $A_1^*$ to $A_2^*$ defined by $f^*(a) = f(a)$ for $a \in A_1$ and for $\pair UV \in A'$
  we put $f^*(\pair UV) = \pair{f(U)}{f(V)}$.
  Then it is easy to show that $f^*$ is a monotone map from $\pair{A_1^*}{\Boxed{\le_1^*}}$ to $\pair{A_2^*}{\Boxed{\le_2^*}}$.
  Moreover, if $f$ is an embedding, then so is $f^*$.
\end{EX}

\begin{EX}
  \textit{A \Katetov\ functor on the category of tournaments and embeddings.}
 Recall that a tournament is a digraph $\pair V E$ such that for every $x,y \in V$ exactly one of the possibilities holds: either $x=y$ or $(x,y) \in E$ or $(y,x) \in E$.

    For a finite set $A$ and a positive integer $n$ let $A^{\le n}$ be the set of all sequences
  $\langle a_1, \ldots, a_k\rangle$ of elements of $A$ where $k \in \{0, 1, \ldots, n\}$.
  In case of $k = 0$ we actually have the empty sequence $\langle\rangle$, as we will be careful to
  distinguish the 1-element sequence $\langle a \rangle$ from $a \in A$.
  For a sequence $s \in A^{\le n}$ let $|s|$ denote the length of~$s$.
  For a tournament $T = \pair V E$, where $E \subseteq V^2$, let $n = |V|$ and let $T^{\le n}$
  be the tournament whose set of vertices is $V^{\le n}$ and whose set of edges is defined \emph{lexicographically}
  as follows:
  \begin{itemize}
  \item
    if $s$ and $t$ are sequences such that $|s| < |t|$, put $s \to t$ in $T^{\le n}$;
  \item
    if $s = \langle s_1, \ldots, s_k\rangle$ and $t = \langle t_1, \ldots, t_k\rangle$ are distinct sequences of the same length,
    find the smallest $i$ such that $s_i \ne t_i$ and then put $s \to t$ in $T^{\le n}$ if and only if $s_i \to t_i$ in $T$.
  \end{itemize}
  For a tournament $T = \pair{V}{E}$ put $K(T) = \pair{V^*}{E^*}$ where
  \begin{align*}
    V^* = V &\mathstrut\union V^{\le n},\\
    E^* = E &\mathstrut\union E(T^{\le n}) \union \{ \pair{v}{s} : v \in V, s \in V^{\le n}, v \text{ appears as an entry in } s\}\\
            &\mathstrut\union \{ \pair{s}{v} : v \in V, s \in V^{\le n}, v \text{ does not appear as an entry in } s\}.
  \end{align*}
  Then it is easy to see that $\pair{V^*}{E^*}$ is a tournament realizing all one-point extensions of $\pair V E$.
  For an embedding $f : \pair{V_1}{E_1} \to \pair{V_2}{E_2}$ let $f^* = K(f)$ be the mapping from
  $V_1^*$ to $V_2^*$ defined by $f^*(v) = f(v)$ for $v \in V_1$ and for $\langle s_1, \ldots, s_k\rangle \in V_1^{\le n}$
  we put $f^*(\langle s_1, \ldots, s_k\rangle) = \langle f(s_1), \ldots, f(s_k)\rangle$.
  Then it is easy to show that $f^*$ is an embedding from $\pair{V_1^*}{E_1^*}$ to $\pair{V_2^*}{E_2^*}$.
  Finally, $K$ is a \Katetov\ functor which is witnessed by the obvious natural transformation mapping $T = \pair V E$ to its copy in $K(T)$.
\end{EX}

\begin{EX}\label{kf.ex.BA}
  \textit{A \Katetov\ functor on the category of all Boolean algebras.}
  For a finite set $A$ let $B(A)$ denote the finite Boolean algebra whose set of atoms is~$A$. For a finite Boolean algebra
  $B(A)$ put $K(B(A)) = B(\{0,1\} \times A)$ and let $\eta_{B(A)} : B(A) \hookrightarrow B(\{0,1\} \times A)$ be the
  unique homomorphism which takes $a \in A$ to $\pair 0a \lor \pair 1a \in B(\{0,1\} \times A)$. Clearly, $\eta_{B(A)}$
  is an embedding. Let us define $K$ on homomorphisms between finite Boolean algebras as follows. Let
  $f : B(A) \to B(A')$ be a homomorphism and assume that for $a \in A$ we have $f(a) = \bigvee S(a)$
  for some $S(a) \subseteq A'$, with the convention that $\bigvee \0 = 0$. Then for $i \in \{0, 1\}$ let
  $K(f)(\pair ia) = \bigvee (\{i\} \times S(a))$. This turns $K$ into a functor from the category of finite
  Boolean algebras into itself which preserves embeddings and such that $\eta : \ID \to K$ is a natural transformation.

  Let us show that $K$ is indeed a \Katetov\ functor. Let $j : B(A) \dothookrightarrow B(A')$.
  Then $B(A') = \langle B(A) \union \{x\}\rangle$ since $B(A')$ is a one-point extension of $B(A)$, so
  $A' = \left(\UNION_{a \in A} \{x \land j(a), \overline x \land j(a)\}\right) \setminus \{0\}$.
  Let $g : B(A') \hookrightarrow K(B(A))$ be the embedding defined on the atoms of $B(A')$ as follows:
  \begin{itemize}
  \item
    if $x \land j(a) = j(a)$ (and consequently $\overline x \land j(a) = 0$)
    or $x \land j(a) = 0$ (and consequently $\overline x \land j(a) = j(a)$)
    let $g$ take $j(a)$ to $\pair 0a \lor \pair 1a$,
  \item
    if $x \land a \ne 0$ and $\overline x \land a \ne 0$ let
    $g$ take $x \land a$ to $\pair 0a$ and $\overline x \land a$ to $\pair 1a$,
  \end{itemize}
  and which extends to the rest of $B(A')$ in the obvious way. Then it is easy to see that $g \circ j = \eta_{B(A)}$.
\end{EX}

\subsection{Sufficient conditions for the existence of \Katetov\ functors}


Let $\Delta$ be a purely relational language,
let $A$ be a $\Delta$-structure, and let $B_1, B_2$ be $\Delta$-structures such that
$A$ is a substructure of both of them and $A = B_1 \cap B_2$.
The \emph{free amalgam} of the $B_1, B_2$ over $A$ is the $\Delta$-structure $C$ with universe
$B_1 \cup B_2$ such that both $B_1$, $B_2$ are substructures of $C$ and for every $R \in \Delta$ we have that
$R^C = R^{B_1} \cup R^{B_2}$ (in other words,
no tuple which meets $B_1 \setminus A$ and $B_2 \setminus A$ satisfies any relation symbol in~$\Delta$).
Following \cite{bilge-melleray}, we say that $\catA$ has the \emph{free amalgamation property} if every triple $A$, $B_1$, $B_2$ as above has the free amalgam in $\catA$.
The next result is implicit in \cite{bilge-melleray} (see Definition 3.7 in \cite{bilge-melleray} and the comment that follows).

\begin{THM}[implicit in \cite{bilge-melleray}]
  If $\catA$ has free amalgamations then a \Katetov\ functor $K : \catA \to \catC$ exists.
\end{THM}

The following theorem is a strengthening of
this as well as of
the main result of~\cite{dolinka-masulovic}.
We say that $\catA$ \emph{has one-point extension pushouts}
[resp. \emph{mixed pushouts}]
in $\catC$ if for every
morphism $f : A_0 \to A_1$ in $\catA$ and a one-point extension
[resp. embedding]
$g : A_0 \dothookrightarrow A_2$ in $\catA$
there exists a $B \in \Ob(\catA)$,
an embedding
$p : A_1 \hookrightarrow B$ and
a morphism $q : A_2 \to B$ such that $p \circ f = q \circ g$ and this commuting square
is a pushout square in the category
$\catCh$ of \emph{all} homomorphisms between $\catC$-objects.
$$
  \XYMATRIX{
    A_0 \dotaremb[r]_g \ar[d]_{f} & A_2 \ar[d]^q\\
    A_1  \aremb[r]_p & B
  }
$$
Note that free amalgamations are particular examples of pushouts.
Note also that typical categories of models with embeddings rarely have pushouts.
Namely, recall that a pair of morphisms $(p,q)$ \emph{provides the pushout} of $(f,g)$ if $p \cmp f = q \cmp g$ and for every other pair $(p',q')$ satisfying $p' \cmp f = q' \cmp g$ there exists a unique morphism $h$ such that $h \cmp p = p'$ and $h \cmp q = q'$.
Now, if $f=g$ and $(p',q')$ consists of identities then clearly $h$ cannot be an embedding.
That is why, in the definition above, we have to consider pushouts in the category of all homomorphisms.

\begin{LEM}\label{LmgEWRoi}
Suppose
$$
  \XYMATRIX{
    A_0 \aremb[r]_g \ar[d]_{f} & A_2 \ar[d]^q\\
    A_1  \aremb[r]_p & B
  }
$$
is a pushout square in the category $\catCh$ of all homomorphisms.
If $g$ is a one-point extension then so is $p$.
\end{LEM}
\begin{PROOF}
Let $B_1$ be the substructure of $B$ generated by $p[A_1] \cup \{q(s)\}$, where $s \in A_2$ is such that $g[A_0] \cup \{s\}$ generates $A_2$.
Notice that $q[A_2] \subseteq B_1$.
In other words, the square
$$
  \XYMATRIX{
    A_0 \aremb[r]_g \ar[d]_{f} & A_2 \ar[d]^{q_1}\\
    A_1  \aremb[r]_{p_1} & B_1
  }
$$
is commutative, where $p_1$ and $q_1$ denote the same mappings as $p$ and $q$, respectively.
By the universality of a pushout, there is a unique homomorphism $\map h B {B_1}$ such that $h \cmp p = p_1$ and $h \cmp q = q_1$.
Let $h_1$ the composition of $h$ with the inclusion $B_1 \subseteq B$.
Again by the universality of a pushout, $\map {h_1}B B$ is the unique homomorphism satisfying $h_1 \cmp p = p$ and $h_1 \cmp q = q$.
It follows that $h_1 = \id_B$ and hence $B_1 = B$.
This completes the proof.
\end{PROOF}

It turns out that both variants of the definition above are equivalent.
In practice however, it is usually easier to verify the existence of pushouts for one-point extensions.

\begin{PROP}
The following properties are equivalent:
\begin{enumerate}
\item[{\rm(a)}] $\catA$ has the one-point extension pushouts in $\catC$.
\item[{\rm(b)}] $\catA$ has mixed pushouts in $\catC$.
\end{enumerate}
\end{PROP}
\begin{PROOF}
Only implication {\rm(a)}$\implies${\rm(b)} requires a proof.
Fix $f : A_0 \to A_1$ and $g : A_0 \hookrightarrow A_2$ as above and assume that $g = g_n \cmp \dots \cmp g_1$ is the composition of $n$ one-point extensions $g_i : E_i \dothookrightarrow E_{i+1}$, where $E_1=A_0$, $E_{n+1}=A_2$.
By Lemma~\ref{LmgEWRoi} we have the following sequence of pushout squares in $\catCh$.
$$
\XYMATRIX{
A_0 \dotaremb[r]_{g_1} \ar[d]_f & E_2 \dotaremb[r]_{g_2} \ar[d]^{q_1} & E_3 \dotaremb[r] \ar[d]^{q_2} & \cdots \cdots \dotaremb[r] & E_n \dotaremb[r]_{g_n} \ar[d]^{q_{n-1}} & A_2 \ar[d]^{q_n} \\
A_1 \dotaremb[r]_{p_1} & B_2 \dotaremb[r]_{p_2} & B_3 \dotaremb[r] & \cdots \cdots \dotaremb[r] & B_n \dotaremb[r]_{p_n} & B
}
$$
Clearly, the composition of all these squares is a pushout in $\catCh$.
\end{PROOF}

\begin{THM}\label{kf.thm.OEP}
  If $\catA$ has has one-point extension pushouts in $\catC$ then a \Katetov\ functor $K : \catA \to \catC$ exists.
\end{THM}
\begin{PROOF}
  Let us first show that every countable source $(A \dothookrightarrow B_n)_{n \in \NN}$ has a pushout in $\catC$,
  where $A, B_1, B_2, \ldots \in \Ob(\catA)$. Let $e_n : A \dothookrightarrow B_n$ be the embeddings in this source.
  Let $P_2 \in \Ob(\catA)$ together with the embeddings $f_2 : B_1 \hookrightarrow P_2$ and $g_2 : B_2 \hookrightarrow P_2$
  be the pushout of $e_1$ and $e_2$. Next,
  let $P_3 \in \Ob(\catA)$ together with the embeddings $f_3 : P_2 \hookrightarrow P_3$ and $g_3 : B_3 \hookrightarrow P_3$
  be the pushout of $f_2 \circ e_1$ and $e_3$. Then,
  let $P_4 \in \Ob(\catA)$ together with the embeddings $f_4 : P_3 \hookrightarrow P_4$ and $g_4 : B_4 \hookrightarrow P_4$
  be the pushout of $f_3 \circ f_2 \circ e_1$ and $e_4$, and so on:
  $$
    \XYMATRIX{
      B_1 \aremb[r]^{f_2} & P_2 \aremb[r]^{f_3} & P_3 \aremb[r]^{f_4} & P_4 \aremb[r]^{f_5} & \cdots\\
      A \dotaremb[u]_{e_1} \dotaremb[r]_{e_2} \dotaremb@/_6mm/[rr]_{e_3} \dotaremb@/_12mm/[rrr]_{e_4}
      & B_2 \ar[u]_{g_2} & B_3 \ar[u]_{g_3} & B_4 \ar[u]_{g_4}
    }
  $$
  Let $P \in \Ob(\catC)$ be the colimit of the chain $B_1 \hookrightarrow P_2 \hookrightarrow P_3 \hookrightarrow
  P_4 \hookrightarrow \cdots$. It is easy to show that $P$ is the pushout of the
  source $(A \dothookrightarrow B_n)_{n \in \NN}$.

  Let us now construct the \Katetov\ functor as the pushout of all the one-point extensions of an object in $\catA$.
  More precisely, for every $A \in \Ob(\catA)$ let us fix embeddings $e_n : A \dothookrightarrow B_n$, where
  $B_1, B_2, \ldots$ is the list of all the one-point extensions of $A$, where every isomorphism type is taken exactly once to keep
  the list countable. Define $K(A)$ to be the pushout of the source $(e_n : A \dothookrightarrow B_n)_{n}$. This is how $K$
  acts on objects.

  Let us show how $K$ acts on morphisms. Take any morphism $h : A \to A'$ in $\catA$.
  Let $(e_i : A \dothookrightarrow B_i)_{i \in I}$ be the source consisting of all the one-point extensions of $A$
  (with every isomorphism type is taken exactly once), and let
  let $(e'_j : A' \dothookrightarrow B'_j)_{j \in J}$ be the source consisting of all the one-point extensions of $A'$
  (with every isomorphism type is taken exactly once). By the assumption, for every $i \in I$ there exists an $m(i) \in J$
  and a morphism $h_i : B_i \to B'_{m(i)}$ such that the following is a pushout square in $\catC$:
  $$
    \XYMATRIX{
      A \dotaremb[r]_{e_i} \ar[d]_{h} & B_i \ar[d]^-{h_i}\\
      A'  \dotaremb[r]_-{e'_{m(i)}} & B'_{m(i)}
    }
  $$
  Now, $K(A')$ is a pushout of the source $(e'_j : A' \dothookrightarrow B'_j)_{j \in J}$ so let us denote the canonical
  embeddings $B'_j \hookrightarrow K(A')$ by $\iota'_j$, $j \in J$. Therefore,
  $(\iota'_{m(i)} \circ h_i : B_i \to K(A'))_{i \in I}$ is a compatible cone over the
  source $(e_i : A \dothookrightarrow B_i)_{i \in I}$, so there is a unique mediating morphism $\tilde h : K(A) \to K(A')$.
  Then we put $K(h) = \tilde h$.
\end{PROOF}

Note that the category of graphs and homomorphisms has pushouts, while the category of $K_n$-free graphs has pushouts of embeddings only.
On the other hand, categories like tournaments or linear orderings do not have pushouts, even when considering all homomorphisms.

\section{\Katetov\ construction}\label{SectKatConstrs}

\begin{DEF}
  Let $K : \catC \to \catC$ be a \Katetov\ functor.
  A \emph{\Katetov\ construction} is a chain of the form:
  $$
    \XYMATRIX{
      C \aremb[r]^-{\eta_C} & K(C) \aremb[r]^-{\eta_{K(C)}} & K^2(C) \aremb[r]^-{\eta_{K^2(C)}} & K^3(C) \aremb[r] & \cdots
    }
  $$
  where $C \in \Ob(\catC)$. We denote the colimit of this chain by $K^\omega(C)$.
  An object $L \in \Ob(\catC)$ \emph{can be obtained by the \Katetov\ construction starting from $C$} if
  $L = K^\omega(C)$. We say that $L$ \emph{can be obtained by the \Katetov\ construction} if
  $L = K^\omega(C)$ for some $C \in \Ob(\catC)$.
\end{DEF}

Note that $K^\omega$ is actually a functor from $\catC$ into $\catC$. Namely, for a morphism $f : A \to B$
let $K^\omega(f)$ be the unique morphism $K^\omega(A) \to K^\omega(B)$ from the colimit of the \Katetov\ construction
starting from $A$ to the competitive compatible cone with the tip at $K^\omega(B)$ and
morphisms $(\hbox{$\hookrightarrow$} \circ K^n(f))_{n \in \NN}$:
$$
  \XYMATRIX{
    & & & & K^\omega(A) \ar@{.>}[ddd]^{K^\omega(f)}\\
    A \aremb[r]_-{\eta_A} \aremb@/^6mm/[rrrru]^-{\eta^\omega_A} \ar[d]_f & K(A) \aremb[r]_{\eta_{K(A)}} \aremb@/^2mm/[rrru] \ar[d]_{K(f)} & K^2(A) \aremb[r]_-{\eta_{K^2(A)}} \aremb[rru] \ar[d]_{K^2(f)} & \cdots &  \\
    B \aremb[r]^-{\eta_B} \aremb@/_6mm/[rrrrd]_-{\eta^\omega_B} & K(B) \aremb[r]^{\eta_{K(B)}} \aremb@/_2mm/[rrrd] & K^2(B) \aremb[r]^-{\eta_{K^2(B)}} \aremb[rrd] & \cdots &  \\
    & & & & K^\omega(B)
  }
$$
It is clear that $K^\omega$ preserves embeddings (the colimit of embeddings is an embedding).
Moreover, the canonical embeddings $\eta^\omega_A : A \hookrightarrow K^\omega(A)$ constitute a natural transformation
$\eta^\omega : \ID \to K^\omega$. Thus, we have:

\begin{THM}
  $K^\omega : \catC \to \catC$ is a \Katetov\ functor.
\end{THM}

Recall that a countable structure $L$ is \emph{ultrahomogeneous} if every isomorphism between two finitely generated
substructures of $L$ extends to an automorphism of $L$. More precisely, $L$ is ultrahomogeneous if
for all $A, B \in \age(L)$, embeddings $j_A : A \hookrightarrow L$ and $j_B : B \hookrightarrow L$, and
for every isomorphism $f : A \to B$ there is an automorphism $f^*$ of $L$ such that $j_B \circ f = f^* \circ j_A$.
$$
  \XYMATRIX{
    A \ar[d]_f \aremb[r]^{j_A} & L \ar[d]^{f^*} \\
    B \aremb[r]^{j_B} & L
  }
$$
One of the crucial points of the classical \Fraisse\ theory is the fact that every ultrahomogeneous structure is the \Fraisse\ limit of its age, and every \Fraisse\ limit is ultrahomogeneous.

Analogously, we say that a countable structure $L$ is \emph{$\catC$-mor\-phism-homo\-ge\-neous},
if every $\catC$-morphism between two finitely generated substructures of $L$ extends to a $\catC$-endomorphism of $L$.
More precisely, $L$ is $\catC$-morphism-ho\-mo\-ge\-neous if
for all $A, B \in \age(L)$, embeddings $j_A : A \hookrightarrow L$ and $j_B : B \hookrightarrow L$, and
for every $\catC$-morphism $f : A \to B$ there is a $\catC$-endomorphism $f^*$ of $L$ such that $j_B \circ f = f^* \circ j_A$.
In particular, if $\catC$ is the category of all countable $\Delta$-structures with all homomorphisms between them,
instead of saying that $L$ is $\catC$-morphism-homogeneous, we say that $L$ is \emph{homomorphism-homogeneous}.
The study of homomorphism-homogeneity was initiated by Cameron \& Ne\v set\v ril~\cite{CamNes}.

The first part of the next result can be viewed as an analogy to Banach's contraction principle: iterating a \Katetov\ functor, starting from an arbitrary object, one always ``tends" to the \Fraisse\ limit, which can be regarded as a ``fixed point" of the \Katetov\ functor.

\begin{THM}\label{kf.thm.Kat-AP}
  If there exists a \Katetov\ functor $K : \catA \to \catC$, then $\catA$ is an amalgamation class,
  it has a \Fraisse\ limit $L$ in $\catC$, and $L$ can be obtained by the \Katetov\ construction starting from an
  arbitrary $C \in \Ob(\catC)$. Moreover, $L$ is $\catC$-morphism-homogeneous.
\end{THM}
\begin{PROOF}
  Take any $C \in \Ob(\catC)$, let
  \begin{equation}\label{kf.eq.1}
    \XYMATRIX{
      C \aremb[r]^-{\eta_C} & K(C) \aremb[r]^{\eta_{K(C)}} & K^2(C) \aremb[r]^{\eta_{K^2(C)}} & K^3(C) \aremb[r] & \cdots
    }
  \end{equation}
  be the \Katetov\ construction starting from $C$, and let $L \in \Ob(\catC)$ be the colimit of this chain.
  Let $\iota_n : K^n(C) \hookrightarrow L$ be the canonical embeddings of the colimit diagram.

  Let us first show that $\age(L) = \Ob(\catA)$.
  Lemma~\ref{kf.lem.age} yields $\age(L) \subseteq \Ob(\catA)$, so let us show that $\Ob(\catA) \subseteq \age(L)$.
  Take any $B \in \Ob(\catA)$ and let $A_1 \hookrightarrow A_2 \hookrightarrow \cdots$ be a chain whose colimit is~$C$.
  Since $\catA$ has (JEP) there is a $D \in \Ob(\catA)$ such that $A_1 \hookrightarrow D \hookleftarrow B$.
  Lemma~\ref{kf.lem.n-chain} then ensures that there is an $n \in \NN$ such that $D \hookrightarrow K^n(A_1)$.
  On the other hand, $A_1 \hookrightarrow C$ implies $K^n(A_1) \hookrightarrow K^n(C)$. Therefore,
  $B \hookrightarrow D \hookrightarrow K^n(A_1) \hookrightarrow K^n(C) \hookrightarrow L$, so $B \in \age(L)$.
  This completes the proof that $\age(L) = \Ob(\catA)$.

  Next, let us show that $L$ \emph{realizes all one-point extensions}, that is, let us show that for all
  $A, B \in \Ob(\catA)$ such that $A \dothookrightarrow B$ and every embedding $f : A \hookrightarrow L$
  there is an embedding $g : B \hookrightarrow L$ such that:
    \begin{equation}\label{kf.eq.2}
      \XYMATRIX{
        A \aremb[r]^{f} \dotaremb[d] & L \\
        B \aremb[ur]_{g}
      }
    \end{equation}
  Take any $A, B \in \Ob(\catA)$ such that $A \dothookrightarrow B$ and let $f : A \hookrightarrow L$ be
  an arbitrary embedding. By Lemma~\ref{kf.lem.factoring-colim} there is an $n \in \NN$ and an
  embedding $h : A \hookrightarrow K^n(C)$ such that $f \circ h = \iota_n$.
  Note that the following diagram commutes:
  $$
    \XYMATRIX{
      A \dotaremb[d] \aremb[dr]^{\eta_A} \aremb[r]^-h & K^n(C) \aremb[dr]^{\eta_{K^n(C)}} \aremb[r]^-{\iota_n} & L \\
      B \aremb[r]^-j & K(A) \aremb[r]^-{K(h)} & K^{n+1}(C) \aremb[u]_{\iota_{n+1}}
    }
  $$
  (the triangle on the left commutes due to the definition of the \Katetov\ functor,
  the parallelogram in the middle commutes because $\eta$ is a natural transformation, while
  the triangle on the right commutes as part of the colimit diagram for the chain $(\ref{kf.eq.1})$).
  Let $g = \iota_{n+1} \circ K(h) \circ j$. Having in mind that $f = \iota_n \circ h$,
  from the last commuting diagram we immediately get that the diagram $(\ref{kf.eq.2})$ commutes for this particular
  choice of~$g$.

  Therefore, $L$ realizes all one-point extensions, so $L$ is an ultrahomogeneous countable structure
  whose age is $\Ob(\catA)$. Consequently, $L$ is the \Fraisse\ limit of $\Ob(\catA)$, whence we easily
  conclude that $\catA$ is an amalgamation class. Moreover, the \Fraisse\ limit of $\catA$
  can be obtained by the \Katetov\ construction starting from an
  arbitrary $C \in \Ob(\catC)$.

  Finally, let us show that $L$ is $\catC$-morphism-homogeneous. Take any $A, B \in \age(L)$, fix embeddings
  $j_A : A \hookrightarrow L$ and $j_B : B \hookrightarrow L$, and let $f : A \to B$ be an arbitrary
  morphism. Then
  \begin{equation}\label{kf.eq.3}
    \XYMATRIX{
      A \ar[rr]^f \aremb[d]_{\eta^\omega_A} & & B \aremb[d]^{\eta^\omega_B} \\
      K^\omega(A) \ar[rr]_{K^\omega(f)}     & & K^\omega(B)
    }
  \end{equation}
  Having in mind that $K^\omega(A)$ and $K^\omega(B)$ are colimits of \Katetov\ constructions starting from $A$ and $B$,
  respectively, we conclude that both $K^\omega(A)$ and $K^\omega(B)$ are isomorphic to $L$. Since $L$ is
  ultrahomogeneous, there exist isomorphisms $s : K^\omega(A) \to L$ and $t : K^\omega(B) \to L$ such that
  \begin{equation}\label{kf.eq.4}
    \XYMATRIX{
      & A \aremb[dl]_{\eta^\omega_A} \aremb[d]^{j_A} & & B \aremb[dr]^{\eta^\omega_B} \aremb[d]_{j_B} & \\
      K^\omega(A) \ar[r]_-s & L & & L & K^\omega(B) \ar[l]^-t
    }
  \end{equation}
  Putting diagrams $(\ref{kf.eq.3})$ and $(\ref{kf.eq.4})$ together we obtain
  $$
    \XYMATRIX{
      & A \aremb[dl]_{\eta^\omega_A} \aremb[d]^{j_A} \ar[rr]^f & & B \aremb[dr]^{\eta^\omega_B} \aremb[d]_{j_B} & \\
      K^\omega(A) \ar[r]^-s \ar@/_5mm/[rrrr]_{K^\omega(f)} & L \ar@{.>}[rr]^{f^*} & & L & K^\omega(B) \ar[l]_-t
    }
  $$
  whence follows that $f^* = t \circ K^\omega(f) \circ s^{-1}$ is a $\catC$-endomorphism of $L$
  which extends~$f$. So, $L$ is $\catC$-morphism-homogeneous.
\end{PROOF}

Consequently, if the \Katetov\ functor is defined on a category of countable $\Delta$-structures and all homomorphisms
between $\Delta$-structures, the \Fraisse\ limit of $\catA$ is both ultrahomogeneous and homomorphism-homogeneous.

\begin{EX}
  Let $n \ge 3$ be an integer,
  let $\catC_n$ be the category of all countable $K_n$-free graphs together with all graph homomorphisms,
  and let $\catA_n$ be the full subcategory of $\catC_n$ spanned by all finite $K_n$-free graphs.
  Then there does not exist a \Katetov\ functor $K : \catA_n \to \catC_n$, for if there were one, the Henson graph $H_n$ --
  the \Fraisse\ limit of $\catA_n$ -- would be homomorphism-homogeneous, and we know this is not the case.

  (Proof. Since $H_n$ is universal for all finite $K_n$-free graphs,
  it embeds both $K_{n-1}$ and the star $S_n$, which is the graph where one vertex is adjacent to $n-1$ independent vertices.
  Let $f$ be a partial homomorphism of $H_n$ which maps the $n-1$ independent vertices of the star $S_n$ onto the vertices of $K_{n-1}$.
  If $H_n$ were homomorphism-homogeneous, $f$ would extend to an endomorphism $f^*$ of $H_n$, so $f^*$ applied to the center of the star $S_n$
  would produce a vertex adjacent to each of the vertices of $K_{n-1}$ inducing thus a $K_n$ in $H_n$, which is not possible.)

  Note however that there exists a \Katetov\ functor from the category $\catA'_n$ of all finite $K_n$-free graphs together with all graph
  embeddings to the category $\catC'_n$ of all countable $K_n$-free graphs together with all graph embeddings (see Example~\ref{kf.ex.Kn-free}).
\end{EX}

\subsection{Characterizations of the existence of a \Katetov\ functor}

The following theorem gives a necessary and sufficient condition for a \Katetov\ functor to exist. It depends on a condition
that resembles the Herwig-Lascar-Solecki property (see \cite{herwig-lascar, solecki}).

\begin{DEF}
  A \emph{partial morphism} of $C \in \Ob(\catC)$ is a triple $\langle A, f, B\rangle$ where
  $A, B \le C$ are finitely generated and $f : A \to B$ is a $\catC$-morphism.
  We say that $C \in \Ob(\catC)$ has the \emph{morphism extension property in~$\catC$} if for any choice
  $f_1, f_2, \ldots$ of partial morphisms of $C$ there exist $D \in \Ob(\catC)$ and $m_1, m_2, \ldots \in \End(D)$
  such that $C$ is a substructure of $D$,
  $m_i$ is an extension of $f_i$ for all $i$, and the following \emph{coherence} conditions are satisfied for all
  $i$, $j$ and $k$:
  \begin{itemize}
  \item
    if $f_i = \langle A, \id_A, A\rangle$ then $m_i = \id_D$,
  \item
    if $f_i$ is an embedding, then so is $m_i$, and
  \item
    if $f_i \circ f_j = f_k$ then $m_i \circ m_j = m_k$.
  \end{itemize}
  We say that $\catC$ has the \emph{morphism extension property} if every $C \in \Ob(\catC)$ has the
  morphism extension property in~$\catC$.
\end{DEF}

\begin{THM}\label{kf.thm.KF-MEP}
  The following are equivalent:
  \begin{enumerate}
  \item
    there exists a \Katetov\ functor $K : \catA \to \catC$;
  \item
    $\catA$ has (AP) and $\catC$ has the morphism extension property;
  \item
    $\catA$ has (AP) and the \Fraisse\ limit of $\catA$ has the morphism extension property in~$\catC$.
  \end{enumerate}
\end{THM}
\begin{PROOF}
  $(1) \Rightarrow (2)$:
  From Theorem~\ref{kf.thm.Kat-AP} we know that $\catA$ is an amalgamation class,
  it has a \Fraisse\ limit $L$ in $\catC$, and $L$ can be obtained by the \Katetov\ construction starting from an
  arbitrary $C \in \Ob(\catC)$.
  Now, take any $C \in \Ob(\catC)$ and let us show that $C$ has the morphism extension property in~$\catC$.
  Since $L$ is universal for $\Ob(\catC)$, without loss of generality we can assume that $C \le L$. For every
  finitely generated $A \le C$ fix an isomorphism $j_A : K^\omega(A) \to L$ such that
  $$
    \XYMATRIX{
      A \aremb[r]^-{\eta^\omega_A} \aremb[d]_-{\le} & K^\omega(A) \ar[d]^-{j_A} \\
      C \aremb[r]_-{\le} & L
    }
  $$
  (such an isomorphism exists because $L$ is ultrahomogeneous). Now, for any family
  $\langle A_i, f_i, B_i\rangle$, $i \in I$, of partial morphisms of~$C$ it is easy to see that $L$
  together with its endomorphisms $m_i = j_{B_i} \circ K^\omega(f_i) \circ j_{A_i}^{-1}$, $i \in I$,
  is an extension of~$C$ and its partial morphisms $f_i$, $i \in I$:
  $$
    \XYMATRIX{
      A_i \aremb[r]^-{\eta^\omega_{A_i}} \ar[d]_-{f_i} & K^\omega(A_i) \ar[r]^-{j_{A_i}} \ar[d]_-{K^\omega(f_i)} & L \ar[d]^-{m_i} \\
      B_i \aremb[r]_-{\eta^\omega_{B_i}}               & K^\omega(B_i) \ar[r]_-{j_{B_i}}                         & L
    }
  $$
  The coherence requirements are satisfied since $K^\omega$ is a functor which preserves embeddings.

  $(2) \Rightarrow (3)$:
  Trivial.

  $(3) \Rightarrow (1)$:
  Let $L$ be the \Fraisse\ limit of $\catA$.
  For every $A \in \Ob(\catA)$ fix an embedding $j_A : A \hookrightarrow L$. Then every $\catA$-morphism
  $f : A \to B$ induces a partial morphism $p(f) : j_A(A) \to j_B(B)$ of $L$ by $p(f) = j_B \circ f \circ j_A^{-1}$.
  Since $L$ is a countable structure, there are only countably many partial morphisms $p(f)$, say, $p_1$, $p_2$, \ldots.
  By the assumption of $(3)$ there exist $D \in \Ob(\catC)$ and $m_1, m_2, \ldots \in \End(D)$ such that $L$ is a substructure of $D$,
  $m_i$ is an extension of $p_i$ for all $i$, and the coherence conditions are satisfied.
  Let $e : L \le D$ be the inclusion of $L$ into $D$.

  Define a functor $K : \catA \to \catC$ on objects by $K(A) = D$ and on morphisms by $K(f) = m_i$, where $p(f) = p_i$.
  Let us show that $K$ is indeed a functor. First, note that $K(\id_A) = \id_D = \id_{K(A)}$:
  let $p(\id_A) = p_i$; since $p_i = p(\id_A) = \id_{j_A(A)}$ coherence requirements force that $m_i = \id_D$.
  Then, let us show that $K(g \circ f) = K(g) \circ K(f)$, where $f : A \to B$ and $g : B \to C$.
  Let $k$ and $l$ be positive integers such that $p(f) = p_k = j_B \circ f \circ j_A^{-1}$ and
  $p(g) = p_l = j_C \circ g \circ j_B^{-1}$. Let $s$ be an integer such that $p_s = j_C \circ g \circ f \circ j_A^{-1}$.
  Then $p_l \circ p_k = p_s$, so the coherence requirements imply that $m_l \circ m_k = m_s$. Finally,
  $K(g \circ f) = m_s = m_l \circ m_k = K(g) \circ K(f)$.
  The coherence requirements also ensure that $K$ preserves embeddings.

  Let us now show that the set of arrows $\eta_A = e \circ j_A$ constitutes a natural transformation $\eta : \ID \to K$.
  Take any $\catA$-morphism $f : A \to B$. Then $p(f) = p_i = j_B \circ f \circ j_A^{-1}$
  is a partial morphism of $L$ whose extension is $m_i$. This is why
  the following diagram commutes (where the dashed arrow indicates a partial morphism):
  $$
    \XYMATRIX{
      A \aremb@/^6mm/[rr]^-{\eta_A} \ar[d]_-f \aremb[r]_-{j_A} & L \ar@{-->}[d]^-{p_i} \aremb[r]_-e & D \ar[d]^-{m_i = K(f)} \\
      B \aremb@/_6mm/[rr]_-{\eta_B}           \aremb[r]^-{j_B} & L \aremb[r]^-e & D
    }
  $$

  Finally, let us show that $K(A)$ embeds all one-point extensions of $A$. Let $A \dothookrightarrow B$. Then there is
  an $h : B \hookrightarrow L$ such that
  $$
    \XYMATRIX{
      A \dotaremb[d] \aremb[r]^-{j_A} & L \\
      B \aremb[ur]_{h}
    }
  $$
  since $L$ is the \Fraisse\ limit of $\catA$. Therefore,
  $$
    \XYMATRIX{
      A \aremb@/^6mm/[rr]^-{\eta_A} \dotaremb[d] \aremb[r]_-{j_A} & L \aremb[r]_-e & D \rlap{$\mathstrut = K(A)$}\\
      B \aremb[ur]_{h} \aremb@/_3mm/[urr]_{e \circ h} &  &
    }
  $$
  which concludes the proof.
\end{PROOF}

Note that the Henson graph $H_n$, $n \ge 3$, does not have the morphism extension property with respect to all
graph homomorphisms (for otherwise there would be a \Katetov\ functor defined on the category of all finite $K_n$-free
graphs and all graph homomorphisms, and we know that such a functor cannot exist).

\paragraph{Conjecture.} Every \Fraisse\ limit has the morphism extension property with respect to embeddings.

\vspace{5mm}

The following theorem shows that the existence of a \Katetov\ functor for varieties of algebras
understood as categories whose objects are the algebras of the variety and morphisms are embeddings
is equivalent to the amalgamation property for the category of finitely generated algebras of the variety.

\begin{THM}\label{kf.thm.Kat-AP-UA}
  Let $\Delta$ be an algebraic language and
  let $\catV$ be a variety of $\Delta$-algebras understood as a category whose objects are $\Delta$-algebras
  and morphisms are embeddings. Let $\catA$ be the full subcategory of $\catV$
  spanned by all finitely generated algebras in $\catV$ and let $\catC$ be the full subcategory of $\catV$
  spanned by all countably generated algebras in $\catV$. Assume additionally that there are only countably
  many isomorphism types in $\catA$. Then there exists a \Katetov\ functor $K : \catA \to \catC$ if and only if
  $\catA$ is an amalgamation class.
\end{THM}
\begin{PROOF}
  $(\Rightarrow)$
  Immediately from Theorem~\ref{kf.thm.Kat-AP}.

  $(\Leftarrow)$
  Recall that a \emph{partial algebra} consists of a set $A$ and some partial operations on $A$, where a partial operation
  is any partial mapping $A^n \to A$ for some $n$ (see \cite{gratzer-ua} for further reference on partial algebras).
  Clearly, the class of all partial algebras of a fixed type is
  a free amalgamation class because we can simply identify the elements of the
  common subalgebra and leave everything else undefined.

  According to Theorem~\ref{kf.thm.OEP} it suffices to show that $\catA$ has one-point extension pushouts in $\catC$.
  Take any $A_0, A_1, A_2 \in \Ob(\catA)$ such that $A_0$ embeds into $A_1$ and $A_2$ is a one-point extension of $A_0$.
  Without loss of generality we may assume that $A_0 \le A_1$ and $A_0 \le A_2$. Let $G \subseteq A_0$ be a finite set
  which generates $A_0$, choose $x \in A_2 \setminus A_0$ so that $G \union \{x\}$ generates $A_2$
  and let $H$ be a finite set disjoint from $G$ such that $G \union H$ generates $A_1$.
  Let $S = A_1 \oplus_{A_0} A_2$ be the partial algebra which arises as the free amalgam of $A_1$ and $A_2$ over $A_0$
  in the class of all partial $\Delta$-algebras. Since $\catA$ has the amalgamation property, there is a $C \in \Ob(\catA)$
  such that
  $$
    \XYMATRIX{
      A_0 \dotaremb[r]_\le \aremb[d]_\le & A_2 \aremb[d]\\
      A_1  \aremb[r] & C
    }
  $$
  whence follows that $C$ embeds the partial algebra $S$ in the sense of \cite[\S 28]{gratzer-ua}.
  It is a well-known fact (see again \cite[\S 28]{gratzer-ua}) that if
  $P$ is a partial algebra which embeds into some total algebra from~$\catV$ then the free algebra
  $\mathbb{F}_\catV(P)$ exists in $\catV$. Therefore, $\mathbb{F}_\catV(S)$ exists and belongs to $\catV$.
  It is easy to see that $\mathbb{F}_\catV(S)$ is generated by $\{x\} \union G \union H$, so
  $\mathbb{F}_\catV(S)$ is a one-point extension of $A_1$. It clearly embeds $A_2$, so we have that
  $$
    \XYMATRIX{
      A_0 \dotaremb[r]_\le \aremb[d]_\le & A_2 \aremb[d]\\
      A_1  \dotaremb[r] & \mathbb{F}_\catV(S)
    }
  $$
  The universal mapping property, which is the defining property of free algebras, ensures that the above
  commuting square is actually a pushout square in~$\catC$. This completes the proof that
  $\catA$ has one-point extension pushouts in~$\catC$.
\end{PROOF}

\begin{COR}
    A \Katetov\ functor exists for the category of all finite semilattices,
    the category of all finite lattices and for the category of all finite Boolean algebras.
\end{COR}
\begin{PROOF}
    The proof follows immediately from the fact that all the three classes of algebras are well-known
    examples of amalgamation classes.
\end{PROOF}

\subsection{Automorphism groups and endomorphism monoids}

The existence of a \Katetov\ functor enables us to quickly conclude that the automorphism group of the corresponding
\Fraisse\ limit is universal, as is the monoid of $\catC$-endomorphisms.
As an immediate consequence of Theorem~\ref{kf.thm.Kat-AP} we have:

\begin{COR}\label{kf.cor.emb}
  Let $K : \catA \to \catC$ be a \Katetov\ functor and let $L$ be the \Fraisse\ limit of $\catA$
  (which exists by Theorem~\ref{kf.thm.Kat-AP}). Then for every $C \in \Ob(\catC)$:
  \begin{itemize}
  \item
    $\Aut(C) \hookrightarrow \Aut(L)$;
  \item
    $\End_\catC(C) \hookrightarrow \End_\catC(L)$.
  \end{itemize}
\end{COR}
\begin{PROOF}
  Since $K^\omega$ is a functor, we immediately get that $\Aut(C) \hookrightarrow \Aut(K^\omega(C))$ via $f \mapsto K^\omega(f)$ and
  that $\End_\catC(C) \hookrightarrow \End_\catC(K^\omega(C))$ via $f \mapsto K^\omega(f)$.
  But, $K^\omega(C) \cong L$ due to Theorem~\ref{kf.thm.Kat-AP}.
\end{PROOF}

Recall that $\End_\catC(X)$ may be just the set of all embeddings of $X$ into itself, in case other homomorphisms are not in $\catC$.
This is the case, for example, in the class of $K_n$-free graphs, where there is no \Katetov\ functor acting on all homomorphisms.


\begin{COR}
    For the following \Fraisse\ limits $L$ we have that $\Aut(L)$ embeds all permutation groups on a countable set:
    \begin{itemize}
    \item
      the random graph (proved originally in~\cite{henson}),
    \item
      Henson graphs (proved originally in~\cite{henson}),
    \item
      the random digraph,
    \item
      the rational Urysohn space (follows also from\cite{uspenskij}),
    \item
      the random poset,
    \item
      the countable atomless Boolean algebra,
    \item
      the random semilattice,
    \item
      the random lattice,
	\end{itemize}
    For the following \Fraisse\ limits $L$ we have that $\End(L)$ embeds all transformation monoids on a countable set:
    \begin{itemize}
    \item
      the random graph (proved originally in~\cite{bdd}),
    \item
      the random digraph,
    \item
      the rational Urysohn space,
    \item
      the random poset (proved originally in~\cite{d-rp}),
    \item
      the countable atomless Boolean algebra.
    \end{itemize}
\end{COR}
\begin{PROOF}
    Having in mind Corollary~\ref{kf.cor.emb},
    in each case it suffices to show that the corresponding category
    $\calC$ contains a countable structure whose automorphism group embeds $\Sym(\NN)$ and whose
    endomorphism monoid embeds $\NN^\NN$ considered as a transformation monoid.
    For example, in case of the rational Urysohn space it suffices
    to consider the metric space $(\NN, d)$ where $d(m,n) = 1$ for all $m, n \in \NN$, while in the case of the
    random Boolean algebra it suffices to consider the free Boolean algebra on $\aleph_0$ generators.
\end{PROOF}

For some applications it is important to know whether the embeddings mentioned in Corollary~\ref{kf.cor.emb} above are topological embeddings, when $\Aut(X)$ and $\End(X)$ are endowed with the \emph{pointwise topology}, that is, the topology inherited from the power $X^X$, where $X$ carries the discrete topology.
This natural topology makes the composition operation (and the inverse, in case of $\Aut(X)$) continuous.
Note that $\Aut(X) \subseteq \End(X)$ are closed in $X^X$ (not being a homomorphism is witnessed by a finite set).
In case where $X$ is countable, $X^X$ is the well-known \emph{Baire space}, a canonical Polish space,
and $\Aut(X)$ is isomorphic to a closed subgroup of the countable infinite symmetric group $S_\infty$.
The importance of such groups is demonstrated in the pioneering work~\cite{KePeTo} connecting \Fraisse\ theory with general Ramsey theory and topological dynamics.
As we shall see in a moment, \emph{every} \Katetov\ functor embeds hom-sets preserving their pointwise topology.

Given $\catC$-objects $X,Y$, denote by $\catC(X,Y)$ the set of all $\catC$-morphisms from $X$ to $Y$, endowed with the pointwise topology, that is, the topology inherited from the product $X^Y$ with $X$ discrete.
Note that a sequence $f_n \in \catC(X,Y)$ converges to $f \in \catC(X,Y)$ if and only if for every finite set $S \subseteq X$ there is $n_0$ such that $f_n \rest S = f \rest S$ for every $n \ge n_0$.

\begin{PROP}
Let $\map K \catC \catC$ be a \Katetov\ functor.
Then for every $\catC$-objects $X$, $Y$, the mapping
$$\catC(X,Y) \owns f \mapsto K(f) \in \catC(K(X), K(Y))$$
is a topological embedding.
\end{PROP}
\begin{PROOF}
>From the definition of a \Katetov\ functor, we know that the mapping above (which we also denote by $K$) is one-to-one, as $K(f)$ can be viewed as an extension of $f$ (the natural transformation $\eta$ consists of embeddings).
Let $f_n$ be a sequence in $\catC(X,Y)$.
If $K(f_n) \to K(f)$ pointwise, then $f_n \to f$ pointwise, due to the remark above.
Now suppose $f_n \to f$ pointwise and fix $a \in K(X)$.
Choose a finite $S \subseteq X$ such that the structure $A=\gen S{}$ generated by $S$ has the property that $a \in K(A)$, after identifying $K(A)$ with a suitable substructure of $K(X)$ (recall that $K(X)$ is the standard colimit of $K(F)$, where $F$ runs over all finitely generated substructures of $X$).
There is $n_0$ such that $f_n \rest S = f \rest S$ whenever $n \ge n_0$.
Then also $f_n \rest A = f \rest A$ for every $n \ge n_0$.
Hence $K(f_n) \rest K(A) = K(f) \rest K(A)$ whenever $n \ge n_0$, showing that $f_n(a) \to f(a)$ in the discrete topology.
Finally, note that the topology on $\catC(X,Y)$ is always metrizable (and therefore determined by sequences), because $X$ is countably generated and $\catC(X,Y)$ is homeomorphic (via the restriction operator) to a subspace of $Y^G$, where $G$ is a countable set generating $X$.
\end{PROOF}

\begin{COR}
The embeddings appearing in Corollary~\ref{kf.cor.emb} are topological with respect to the pointwise topology.
\end{COR}

\section{Semigroup Bergman property}\label{SecBergmanPrpty}

Following~\cite{maltcev-mitchell-ruskuc}, we say that
a semigroup $S$ is \emph{semigroup Cayley bounded with respect to a generating set $U$} if
$S = U \union U^2 \union \ldots \union U^n$ for some $n \in \NN$. We say that a \emph{semigroup $S$ has the semigroup Bergman
property} if it is semigroup Cayley bounded with respect to every generating set.

A semigroup $S$ has \emph{\Sierpinski\ rank} $n$ if $n$ is the least positive integer such that for any
countable $T \subseteq S$ there exist $s_1, \ldots, s_n \in S$ such that $T \subseteq \langle s_1, \ldots, s_n\rangle$.
If no such $n$ exists, the \Sierpinski\ rank of $S$ is said to be infinite.
A semigroup $S$ is \emph{strongly distorted} if there exists a sequence of natural numbers $l_1, l_2, l_3, \ldots$ and an $N \in \NN$
such that for any sequence $a_1, a_2, a_3, \ldots \in S$ there exist $s_1, \ldots, s_N \in S$ and a sequence of words
$w_1, w_2, w_3, \ldots$ over the alphabet $\{x_1, x_2, \ldots, x_N\}$ such that $|w_n| \le l_n$ and $a_n = w_n(s_1, \ldots, s_N)$
for all~$n$.

\begin{LEM}[\cite{maltcev-mitchell-ruskuc}]\label{kf.lem.SSD}
  If $S$ is a strongly distorted semigroup which is not finitely generated, then $S$ has the Bergman property.
\end{LEM}

It was shown in~\cite{peresse} that $\End(R)$, the endomorphism monoid of the random graph, is
strongly distorted and hence has the semigroup Bergman property since it is not finitely generated.
The idea from~\cite{peresse} was later in~\cite{dolinka} directly generalized to classes of structures with coproducts.
Here, we present a general treatment in the context of classes for which a \Katetov\ functor exists, and where
the (JEP) can be carried out constructively in the sense of the following definition.

\begin{DEF}
  A category $\catC$ \emph{has natural (JEP)} if there exists a covariant functor $F : \catC \times \catC \to \catC$ such that
  \begin{itemize}
  \item
    for all $C, D \in \Ob(\catC)$ there exist embeddings $\lambda_C : C \hookrightarrow F(C, D)$ and
    $\rho_D : D \hookrightarrow F(C, D)$, and
  \item
    for every pair of morphisms $f : C \to C'$ and $g : D \to D'$ the diagram below commutes:
    $$
      \XYMATRIX{
        C  \aremb[r]^-{\lambda_C} \ar[d]_-{f} & F(C, D) \ar[d]^-{F(f,g)} \ar@{<-^{)}}[r]^-{\rho_D} & D \ar[d]_-{g} \\
        C' \aremb[r]_-{\lambda_{C'}}  & F(C', D') \ar@{<-^{)}}[r]_-{\rho_{D'}} & D
      }
    $$
  \end{itemize}
  We also say that $F$ is a \emph{natural (JEP) functor} for $\catC$.

  A category $\catC$ \emph{has retractive natural (JEP)} if $\catC$ has natural (JEP) and the functor $F$ has the
  following additional property: for every $C \in \Ob(\catC)$ there exist morphisms $\rho^*_C, \lambda^*_C : F(C,C) \to C$
  such that $\rho^*_C \circ \rho_C = \id_C = \lambda^*_C \circ \lambda_C$.
\end{DEF}

\begin{REM}\label{kf.rem.cov-F}
  Note that since $F$ is a covariant functor, the following also holds:
  \begin{itemize}
  \item
    $F(\id_C, \id_D) = \id_{F(C, D)}$ for all $C, D \in \Ob(\catC)$,
  \item
    for all $f_1 : B_1 \to C_1$, $f_2 : B_2 \to C_2$, $g_1 : C_1 \to D_1$, $g_2 : C_2 \to D_2$ we have
    $F(g_1 \circ f_1, g_2 \circ f_2) = F(g_1, g_2) \circ F(f_1, f_2)$, and
  \item
    $
      \XYMATRIX{
        A \ar[r]^{f_3} \ar[d]_{f_1} & C \ar[d]^{f_2} \\
        B \ar[r]_{f_4} & D
      }
    $
    and
    $
      \XYMATRIX{
        P \ar[r]^{g_3} \ar[d]_{g_1} & Q \ar[d]^{g_2} \\
        R \ar[r]_{g_4} & S
      }
    $
    implies
    $
      \XYMATRIX{
        F(A,P) \ar[r]^{F(f_3,g_3)} \ar[d]_{F(f_1,g_1)} & F(C,Q) \ar[d]^{F(f_2,g_2)} \\
        F(B,R) \ar[r]_{F(f_4,g_4)} & F(D,S)
      }
    $
  \end{itemize}
\end{REM}

\begin{EX}
  Any category with coproducts (such as the category of graphs, posets, digraphs) has retractive natural (JEP):
  just take $F(C, D)$ to be the coproduct of $C$ and $D$.
\end{EX}

\begin{EX}
  The category of all countable metric spaces with distances in $[0,1]_\QQ = \QQ \sec [0, 1]$ and nonexpansive mappings has retractive natural (JEP):
  take $F(C, D)$ to be the disjoint union of $C$ and $D$ where the distance between any point in $C$ and any
  point in $D$ is~1.

  On the other hand, it is easy to show that the category of all countable metric spaces with distances in $\QQ$ and nonexpansive mappings
  does not have natural (JEP). Suppose, to the contrary,
  that there exists a functor $F$ which realizes the natural (JEP) in this category, let $U$ be the rational Urysohn space and let
  $W = F(U, U)$.
  Let $a_0, b_0 \in U$ be arbitrary but fixed, and let $\delta = d_W(\lambda_U(a_0), \rho_U(b_0))$.
  Take any $a, b \in U$, let $c_a : U \to U : x \mapsto a$ and $c_b : U \to U : x \mapsto b$ be the constant maps and
  put $\Phi = F(c_a, c_b)$. Then $d_W(\lambda_U(a), \rho_U(b)) = d_W(\lambda_U(c_a(a_0)), \rho_U(c_b(b_0))) =
  d_W(\Phi(\lambda_U(a_0)), \Phi(\rho_U(b_0))) \le d_W(\lambda_U(a_0), \rho_U(b_0)) = \delta$, because $\Phi$ is nonexpansive.
  Now, for $a_1, a_2 \in U$ we have $d_U(a_1, a_2) = d_W(\lambda_U(a_1), \lambda_U(a_2)) \le
  d_W(\lambda_U(a_1), \rho_U(b)) + d_W(\lambda_U(a_2), \rho_U(b)) \le 2\delta$. Hence, $\mathrm{diam}(U) \le 2\delta$. Contradiction.
\end{EX}

\begin{EX}
  Let $\Delta$ be the language consisting of function symbols and constant symbols only so that $\Delta$-structures
  are actually $\Delta$-algebras, and assume that $\Delta$ contains a constant symbol~$1$. Then
  the category of $\Delta$-algebras has retractive natural (JEP):
  take $F(C, D)$ to be $C \times D$ where $\lambda_C : c \mapsto \pair{c}{1^D}$,
  $\rho_D : d \mapsto \pair{1^C}{d}$, $\lambda^*_C = \pi_1$ and $\rho^*_D = \pi_2$.
\end{EX}

Our aim in this section is to prove the following theorem:

\begin{THM}\label{kf.thm.Bergman}
  Assume that there exists a \Katetov\ functor $K : \catA \to \catC$ and assume that $\catC$ has retractive natural (JEP).
  Let $L$ be the \Fraisse\ limit of $\catA$ (which exists by Theorem~\ref{kf.thm.Kat-AP}).
 {Assume additionally that there is a retraction $r : K(L) \to L$ such that $r \circ \eta_L = \id_L$.}
  Then $\End_\catC(L)$ is strongly distorted and its \Sierpinski\ rank is at most~5.
  Consequently, if $\End_\catC(L)$ is not finitely generated then it has the Bergman property.
\end{THM}

The proof of the theorem requires some technical prerequisites.
Let us denote the functor which realizes (JEP) in $\catC$ by $(\cdot, \cdot)$ so that
$(C, D)$ denotes its action on objects, and $(f, g)$ its action on morphisms. For objects $C_1, C_2, C_3, \ldots, C_n$ and
morphisms $f, g, f_1, f_2, f_3, \ldots, f_n$ of $\catC$ let
  \begin{align*}
    [C_1, C_2, C_3, \ldots, C_n] &= ((((C_1, C_2), C_3), \ldots), C_n),\\
    [f_1, f_2, f_3, \ldots, f_n] &= ((((f_1, f_2), f_3), \ldots), f_n),\\
    [f, g]_n &= [f, \underbrace{g, \ldots, g}_n], \text { with } [f,g]_0 = f.\\
  \end{align*}
Moreover, let
  \begin{align*}
    L_1 &= L,\\
    L_n &= (L_{n-1}, L) = [\underbrace{L, L, \ldots, L}_n], \text{ for } n \ge 2.
  \end{align*}
Let $C$ denote the colimit of the following chain in $\catC$ with the canonical embeddings denoted by~$\iota_n$:
  $$
    \XYMATRIX{
      L_1 \aremb@/_2mm/[drrr]_(.3){\iota_1} \aremb[r]^-{\lambda_{L_1}} &
      L_2 \aremb[drr]_(.3){\iota_2} \aremb[r]^-{\lambda_{L_2}} &
      L_3 \aremb[dr]_(.3){\iota_3} \aremb[r]^-{\lambda_{L_3}} & \cdots
    \\
      &                                &                                & C
    }
  $$
Let $L$ be the \Fraisse\ limit of $\catA$, which exists by Theorem~\ref{kf.thm.Kat-AP}.
We know that $K^\omega(C) \cong L$, so let us fix an isomorphism
  $$
    \alpha : K^\omega(C) \mathrel{\overset\cong\longrightarrow} L.
  $$

The following diagram commutes because $(\cdot, \cdot)$ is a natural (JEP) functor:
  $$
    \XYMATRIX{
      L_1 \aremb[r]^-{\lambda_{L_1}} \ar[dr]_-{\id_L}
      &
      L_2 \ar[d]_{\lambda^*_L} \aremb[rr]^-{\lambda_{L_2}} &
      &
      L_3 \ar[d]_{[\lambda^*_L, \id_L]_1} \aremb[rr]^-{\lambda_{L_3}} &
      &
      L_4 \ar[d]_{[\lambda^*_L, \id_L]_2} \aremb[r]^-{\lambda_{L_4}} &
      \cdots
    \\
      &
      L_1 \aremb[rr]_-{\lambda_{L_1}} &
      &
      L_2 \aremb[rr]_-{\lambda_{L_2}} &
      &
      L_3 \aremb[r]_-{\lambda_{L_3}} &
      \cdots\\
    }
  $$
so the following diagram also commutes:
  $$
    \XYMATRIX{
      L_1 \aremb[r]^-{\lambda_{L_1}} \ar[dr]_-{\id_L}
      &
      L_2 \ar[d]_{\lambda^*_L} \aremb[rr]^-{\lambda_{L_2}} &
      &
      L_3 \ar[d]_{[\lambda^*_L, \id_L]_1} \aremb[rr]^-{\lambda_{L_3}} &
      &
      L_4 \ar[d]_{[\lambda^*_L, \id_L]_2} \aremb[r]^-{\lambda_{L_4}} &
      \cdots
    \\
      &
      L_1 &
      &
      L_2 \ar[ll]^-{\lambda^*_{L_1}} &
      &
      L_3 \ar[ll]^-{\lambda^*_{L_2}} &
      \ar[l]^-{\lambda^*_{L_3}} \cdots\\
    }
  $$
Therefore, there is a compatible cone with the tip at $L$ and the morphisms $\id_L$, $\lambda^*_L$,
$\lambda^*_{L_1} \circ [\lambda^*_L, \id_L]_1$, $\lambda^*_{L_1} \circ \lambda^*_{L_2} \circ [\lambda^*_L, \id_L]_2$ \ldots\
over the chain $L_1 \hookrightarrow L_2 \hookrightarrow L_3 \hookrightarrow \cdots$.
Since $C$ is a colimit of the chain, there is a unique $\beta : C \to L$ such that
  $$
    \XYMATRIX{
      C \ar@{.>}[rrrrrr]^\beta&
      &
      &
      &
      &
      &
      L &
    \\
      L_1 \aremb[rr]_-{\lambda_{L_1}} \aremb[u]^-{\iota_1} \aremb@/^2mm/[urrrrrr]^(.6){\id_L} &
      &
      L_2 \aremb[rr]_-{\lambda_{L_2}} \aremb[ull]_(0.55){\iota_2} \aremb[urrrr]_(.575){\lambda^*_L} &
      &
      L_3 \aremb[rr]_-{\lambda_{L_3}} \aremb@/_2mm/[ullll]_{\iota_3} \aremb[urr]_(.6){\lambda^*_{L_1} \circ [\lambda^*_L, \id_L]_1} &
      & \cdots &\\
    }
  $$
In particular,
  \begin{equation}\label{kf.eq.beta-iota1}
    \beta \circ \iota_1 = \id_L.
  \end{equation}

As the next step in the construction, note that the following diagram commutes (again due to the fact that $(\cdot, \cdot)$ is a natural (JEP)
functor):
  $$
    \XYMATRIX{
      L_1 \ar[d]_-{\rho_L} \aremb[rr]^-{\lambda_{L_1}} &
      &
      L_2 \ar[d]_{[\rho_L, \id_L]_1} \aremb[rr]^-{\lambda_{L_2}} &
      &
      L_3 \ar[d]_{[\rho_L, \id_L]_2} \aremb[rr]^-{\lambda_{L_3}} &
      &
      L_4 \ar[d]_{[\rho_L, \id_L]_3} \aremb[r]^-{\lambda_{L_4}} &
      \cdots
    \\
      L_2 \ar[d]_-{\rho^*_L} \aremb[rr]_-{\lambda_{L_2}} &
      &
      L_3 \ar[d]_{[\rho^*_L, \id_L]_1} \aremb[rr]_-{\lambda_{L_3}} &
      &
      L_4 \ar[d]_{[\rho^*_L, \id_L]_2} \aremb[rr]_-{\lambda_{L_4}} &
      &
      L_5 \ar[d]_{[\rho^*_L, \id_L]_3} \aremb[r]_-{\lambda_{L_5}} &
      \cdots
    \\
      L_1 \aremb[rr]_-{\lambda_{L_1}} &
      &
      L_2 \aremb[rr]_-{\lambda_{L_2}} &
      &
      L_3 \aremb[rr]_-{\lambda_{L_3}} &
      &
      L_4 \aremb[r]_-{\lambda_{L_4}} &
      \cdots
    }
  $$
Therefore, there is a compatible cone with the tip at $C$ and the morphisms $\iota_2 \circ \rho_L$,
$\iota_3 \circ [\rho_L, \id_L]_1$, $\iota_4 \circ [\rho_L, \id_L]_2$ \ldots\
over the chain $L_1 \hookrightarrow L_2 \hookrightarrow L_3 \hookrightarrow \cdots$.
Since $C$ is a colimit of the chain, there is a unique $\sigma : C \to C$ such that
  $$
    \XYMATRIX{
      C \ar@{.>}[rrrrrr]^\sigma&
      &
      &
      &
      &
      &
      C &
    \\
      L_1 \aremb[rr]_-{\lambda_{L_1}} \aremb[u]^-{\iota_1} \aremb@/^3mm/[urrrrrr]^(.5){\iota_2 \circ \rho_L} &
      &
      L_2 \aremb[rr]_-{\lambda_{L_2}} \aremb[ull]_(0.55){\iota_2} \aremb@/^3mm/[urrrr]_(.575){\iota_3 \circ [\rho_L, \id_L]_1} &
      &
      L_3 \aremb[rr]_-{\lambda_{L_3}} \aremb@/_2mm/[ullll]_(.6){\iota_3} \aremb[urr]_(.6){\iota_4 \circ [\rho_L, \id_L]_2} &
      & \cdots &\\
    }
  $$
or, explicitly,
  $$
    \sigma \circ \iota_n = \iota_{n+1} \circ [\rho_L, \id_L]_{n-1}, \text{ for all } n \ge 1.
  $$
An easy induction on $n$ then suffices to show that
\begin{equation}\label{kf.eq.sigma-n}
  \sigma^n \circ \iota_1 = \iota_{n+1} \circ [\rho_L, \id_L]_{n-1} \circ \ldots \circ [\rho_L, \id_L]_{1} \circ \rho_L,
  \text{ for all } n \ge 1.
\end{equation}
Also, there is a compatible cone with the tip at $C$ and the morphisms $\iota_1 \circ \rho^*_L$,
$\iota_2 \circ [\rho^*_L, \id_L]_1$, $\iota_3 \circ [\rho^*_L, \id_L]_2$ \ldots\
over the chain $L_2 \hookrightarrow L_3 \hookrightarrow L_4 \hookrightarrow \cdots$,
so there is a unique $\tau : C \to C$ such that
  $$
    \XYMATRIX{
      C \ar@{.>}[rrrrrr]^\tau&
      &
      &
      &
      &
      &
      C &
    \\
      L_2 \aremb[rr]_-{\lambda_{L_1}} \aremb[u]^-{\iota_2} \aremb@/^3mm/[urrrrrr]^(.5){\iota_1 \circ \rho^*_L} &
      &
      L_3 \aremb[rr]_-{\lambda_{L_2}} \aremb[ull]_(0.55){\iota_3} \aremb@/^3mm/[urrrr]_(.575){\iota_2 \circ [\rho^*_L, \id_L]_1} &
      &
      L_4 \aremb[rr]_-{\lambda_{L_3}} \aremb@/_2mm/[ullll]_(.6){\iota_4} \aremb[urr]_(.6){\iota_3 \circ [\rho^*_L, \id_L]_2} &
      & \cdots &\\
    }
  $$
or, explicitly,
$$
    \tau \circ \iota_{n+1} = \iota_{n} \circ [\rho^*_L, \id_L]_{n-1}, \text{ for all } n \ge 1.
$$
Another easy induction on $n$ suffices to show that
\begin{equation}\label{kf.eq.tau-n}
    \tau^n \circ \iota_{n+1} = \iota_{1} \circ \rho^*_L \circ [\rho^*_L, \id_L]_{1} \circ \ldots \circ [\rho^*_L, \id_L]_{n-1}, \text{ for all } n \ge 1.
\end{equation}

Let $\overline f = (f_1, f_2, \ldots)$ be a sequence of $\catC$-endomorphisms of $L$. As the final step,
we shall now construct an endomorphism $\phi(\overline f) : C \to C$ which encodes the sequence $\overline f$.
Using once more the fact that $(\cdot, \cdot)$ is a natural (JEP) functor, we immediately get that
the following diagram commutes:
$$
    \XYMATRIX{
      L_1 \ar[d]_-{f_1} \aremb[rr]^-{\lambda_{L_1}} &
      &
      L_2 \ar[d]_{[f_1, f_2]} \aremb[rr]^-{\lambda_{L_2}} &
      &
      L_3 \ar[d]_{[f_1, f_2, f_3]} \aremb[rr]^-{\lambda_{L_3}} &
      &
      L_4 \ar[d]_{[f_1, f_2, f_3, f_4]} \aremb[r]^-{\lambda_{L_4}} &
      \cdots
    \\
      L_1 \aremb[rr]_-{\lambda_{L_1}} &
      &
      L_2 \aremb[rr]_-{\lambda_{L_2}} &
      &
      L_3 \aremb[rr]_-{\lambda_{L_3}} &
      &
      L_4 \aremb[r]_-{\lambda_{L_4}} &
      \cdots
    }
$$
so there is a unique $\phi(\overline f) : C \to C$ such that
  $$
    \XYMATRIX{
      C \ar@{.>}[rrrrrr]^{\phi(\overline f)}&
      &
      &
      &
      &
      &
      C &
    \\
      L_1 \aremb[rr]_-{\lambda_{L_1}} \aremb[u]^-{\iota_1} \aremb@/^3mm/[urrrrrr]^(.5){\iota_1 \circ f_1} &
      &
      L_2 \aremb[rr]_-{\lambda_{L_2}} \aremb[ull]_(0.55){\iota_2} \aremb@/^3mm/[urrrr]_(.575){\iota_2 \circ [f_1, f_2]} &
      &
      L_3 \aremb[rr]_-{\lambda_{L_3}} \aremb@/_2mm/[ullll]_(.6){\iota_3} \aremb[urr]_(.6){\iota_3 \circ [f_1, f_2, f_3]} &
      & \cdots &\\
    }
  $$
or, explicitly,
  $$
    \phi(\overline f) \circ \iota_n = \iota_{n} \circ [f_1, f_2, \ldots, f_n], \text{ for all } n \ge 1.
  $$

\begin{LEM}\label{kf.lem.calc-1}
  \begin{itemize}
  \item[$(a)$]
    $\phi(\overline f) \circ \iota_1 = \iota_1 \circ f_1$;
  \item[$(b)$]
    $\phi(\overline f) \circ \iota_2 \circ \rho_L= \iota_1 \circ \rho_L \circ f_2$;
  \item[$(c)$]
    $\phi(\overline f) \circ \iota_n \circ [\rho_L, \id_L]_{n-2} \circ \ldots \circ [\rho_L, \id_L]_1 \circ \rho_L =
      \iota_n \circ [\rho_L, \id_L]_{n-2} \circ \ldots \circ [\rho_L, \id_L]_1 \circ \rho_L \circ f_n$, for all $n \ge 3$.
  \end{itemize}
\end{LEM}
\begin{PROOF}
  $(a)$ This is immediate from the construction of $\phi(\overline f)$.

  $(b)$ It suffices to note that the diagram below commutes. The square on the left commutes because $(\cdot,\cdot)$ is
  natural, while the square on the right commutes by the construction of $\phi(\overline f)$.
  $$
    \XYMATRIX{
      L_1 \ar[d]_{f_2} \aremb[r]^{\rho_L} &
      L_2 \ar[d]_{[f_1, f_2]} \aremb[r]^{\iota_2} &
      C \ar[d]^{\phi(\overline f)}
    \\
      L_1 \aremb[r]^{\rho_L} &
      L_2 \aremb[r]^{\iota_2} &
      C
    }
  $$

  $(c)$ This follows by induction on $n$. Just to illustrate the main ideas (which are straightforward, anyhow)
  we show the case $n = 4$. The following diagram commutes:
  $$
    \XYMATRIX{
      L_1 \ar[d]_-{f_4} \aremb[rr]^-{\rho_L} &
      &
      L_2 \ar[d]_{[f_3, f_4]} \aremb[rr]^-{[\rho_L, \id_L]_1} &
      &
      L_3 \ar[d]_{[f_2, f_3, f_4]} \aremb[rr]^-{[\rho_L, \id_L]_2} &
      &
      L_4 \ar[d]_{[f_1, f_2, f_3, f_4]} \aremb[r]^-{\iota_4} &
      C \ar[d]^\phi
    \\
      L_1 \aremb[rr]_-{\rho_L} &
      &
      L_2 \aremb[rr]_-{[\rho_L, \id_L]_1} &
      &
      L_3 \aremb[rr]_-{[\rho_L, \id_L]_2} &
      &
      L_4 \aremb[r]_-{\iota_4} &
      C
    }
  $$
  The leftmost square commutes because $(\cdot,\cdot)$ is
  natural, while the rightmost square commutes by the construction of $\phi(\overline f)$.
  To see that the second square in this row commutes, just apply the functor $(\cdot, \cdot)$ to the following two
  commutative squares (see Remark~\ref{kf.rem.cov-F}):
  $$
    \XYMATRIX{
      L_1 \aremb[r]^-{\rho_L} \ar[d]_{f_3} & L_2 \ar[d]^{[f_2,f_3]} & & L \ar[d]_{f_4} \ar[r]^{\id_L} & L \ar[d]_{f_4} \\
      L_1 \aremb[r]^-{\rho_L} & L_2 & & L \ar[r]^{\id_L} & L
    }
  $$
  The same argument suffices to show that the third square in the row commutes too.
\end{PROOF}

\begin{LEM}\label{kf.lem.calc-2}
  \begin{itemize}
  \item[$(a)$]
    $\beta \circ \phi(\overline f) \circ \iota_1 = f_1$;
  \item[$(b)$]
    $\beta \circ \tau^n \circ \phi(\overline f) \circ \sigma^n \circ \iota_1  = f_{n+1}$.
  \end{itemize}
\end{LEM}
\begin{PROOF}
  In order to make it easier to follow the calculations we underline the expression that
  is to be reduced in the following step.

\medskip

  $(a)$ $\beta \circ \reducenext{\phi(\overline f) \circ \iota_1} = \reducenext{\beta \circ \iota_1} \circ f_1 = f_1$, by Lemma~\ref{kf.lem.calc-1}
  and $(\ref{kf.eq.beta-iota1})$.

\medskip

  $(b)$ $\beta \circ \tau^n \circ \phi(\overline f) \circ (\sigma^n \circ \iota_1) = $
  \begin{alignat*}{5}
  \text{[by $(\ref{kf.eq.sigma-n})$]}
  &= \beta \circ \tau^n \circ \reducenext{\phi(\overline f) \circ \iota_{n+1} \circ [\rho_L, \id_L]_{n-1} \circ \ldots \circ [\rho_L, \id_L]_{1} \circ \rho_L} \\
  \text{[Lemma~\ref{kf.lem.calc-1}]}
  &= \beta \circ \reducenext{\tau^n \circ \iota_{n+1}} \circ [\rho_L, \id_L]_{n-1} \circ \ldots \circ [\rho_L, \id_L]_{1} \circ \rho_L \circ f_{n+1}\\
  \text{[by $(\ref{kf.eq.tau-n})$]}
  &= \reducenext{\beta \circ \iota_{1}} \circ \rho^*_L \circ [\rho^*_L, \id_L]_{1} \circ \ldots \circ [\rho^*_L, \id_L]_{n-1} \circ\\
  & \hphantom{= \beta \circ \iota_{1} \circ} \circ [\rho_L, \id_L]_{n-1} \circ \ldots \circ [\rho_L, \id_L]_{1} \circ \rho_L \circ f_{n+1}\\
  \text{[by $(\ref{kf.eq.beta-iota1})$]}
  &= \rho^*_L \circ [\rho^*_L, \id_L]_{1} \circ \ldots \circ \reducenext{[\rho^*_L, \id_L]_{n-1} \circ [\rho_L, \id_L]_{n-1}} \circ\\
  & \hphantom{= \beta \circ \iota_{1} \circ} \circ [\rho_L, \id_L]_{n-2} \circ \ldots \circ [\rho_L, \id_L]_{1} \circ \rho_L \circ f_{n+1}\\
  &= \ldots = f_{n+1},
  \end{alignat*}
  since $[\rho^*_L, \id_L]_{j} \circ [\rho_L, \id_L]_{j} = \id_L$, for all $j$.
\end{PROOF}

We are now ready to prove Theorem~\ref{kf.thm.Bergman}.

\medskip

\begin{PROOF} \textit{(of Theorem~\ref{kf.thm.Bergman})}
  We are going to show that $\End(K^\omega(C))$, which is isomorphic to $\End(L)$ because $L \cong K^\omega(C)$,
  is strongly distorted and that the \Sierpinski\ rank of $\End(K^\omega(C))$ is at most~5. Take any
  countable sequence $f_1, f_2, \ldots \in \End(K^\omega(C))$, and let us construct $\tilde\alpha$, $\tilde\beta$, $\tilde\sigma$,
  $\tilde\tau$, $\tilde\phi \in \End(K^\omega(C))$ as follows, with the notation introduced above.

  Let $\tilde\alpha = \eta^\omega_C \circ \iota_1 \circ \alpha : K^\omega(C) \to K^\omega(C)$.
  We shall construct $\tilde\beta : K^\omega(C) \to K^\omega(C)$ such that $\tilde\beta \circ \eta^\omega_C = \alpha^{-1} \circ \beta$.
  Since $\eta$ is natural, the diagram on the left below commutes, so by taking $\beta_1 = r \circ K(\beta)$ we have that the
  diagram on the right also commutes:
  $$
    \XYMATRIX{
      C \ar[r]^-\beta \aremb[d]_-{\eta_C} & L \aremb[d]_-{\eta_L} & & C \ar[r]^-\beta \aremb[d]_-{\eta_C} & L\\
      K(C) \ar[r]_-{K(\beta)}             & K(L)                  & & K(C) \ar[ur]_-{\beta_1}
    }
  $$
  Analogously, the following diagrams also commute where $\beta_2 = r \circ K(\beta_1)$:
  $$
    \XYMATRIX{
      K(C) \ar[r]^-{\beta_1} \aremb[d]_-{\eta_{K(C)}} & L \aremb[d]_-{\eta_L} & & K(C) \ar[r]^-{\beta_1} \aremb[d]_-{\eta_{K(C)}} & L\\
      K^2(C) \ar[r]_-{K(\beta_1)}                     & K(L)                  & & K^2(C) \ar[ur]_-{\beta_2}
    }
  $$
  And so on. We get a sequence of morphisms $\beta_n : K^n(C) \to L$ such that
  $$
    \XYMATRIX{
      K^n(C) \ar[r]^-{\beta_n} \aremb[d]_-{\eta_{K^n(C)}} & L\\
      K^{n+1}(C) \ar[ur]_-{\beta_{n+1}}
    }
  $$
  Since $K^\omega(C)$ is the colimit of the chain
  $$
    \XYMATRIX{
      C \aremb[r]^-{\eta_C} & K(C) \aremb[r]^-{\eta_{K(C)}} & K^2(C) \aremb[r]^-{\eta_{K^2(C)}} & K^3(C) \aremb[r] & \cdots
    }
  $$
  there is a unique mediating morphism $\beta_\omega : K^\omega(C) \to L$ such that
  $$
    \XYMATRIX{
      K^\omega(C) \ar@/^8mm/[rr]^-{\beta_\omega}
                  & K^n(C) \ar[r]^-{\beta_n} \aremb[d]_-{\eta_{K^n(C)}} \aremb[l] & L\\
                  & K^{n+1}(C) \ar[ur]_-{\beta_{n+1}} \aremb[ul]
    }
  $$
  In particular, $\beta_\omega \circ \eta^\omega_C = \beta$. Now put $\tilde\beta = \alpha^{-1} \circ \beta_\omega$.

  Finally, let $\tilde\sigma = K^\omega(\sigma)$ and $\tilde\tau = K^\omega(\tau)$,
  let $f^\alpha_n = \alpha \circ f_n \circ \alpha^{-1}$,
  and let $\tilde\phi = K^\omega(\phi(\overline g))$ where $\overline g = (f^\alpha_1, f^\alpha_2, \ldots)$. Then
  \begin{align*}
    \tilde\beta \circ \tilde\phi \circ \tilde\alpha
    &= \tilde\beta \circ \reducenext{K^\omega(\phi(\overline g)) \circ \eta^\omega_C} \circ \iota_1 \circ \alpha\\
  [\text{$\eta^\omega$ is natural}]
    &= \reducenext{\tilde\beta \circ \eta^\omega_C} \circ \phi(\overline g) \circ \iota_1 \circ \alpha\\
  [\text{definition of $\tilde\beta$}]
    &= \alpha^{-1} \circ \reducenext{\beta \circ \phi(\overline g) \circ \iota_1} \circ \alpha\\
  [\text{Lemma~\ref{kf.lem.calc-2}}]
    &= \alpha^{-1} \circ f^\alpha_1 \circ \alpha = f_1,
  \end{align*}
  and
  \begin{align*}
    \tilde\beta \circ \tilde\tau^n \circ \tilde\phi \circ \tilde\sigma^n \circ \tilde\alpha
    &= \tilde\beta \circ \reducenext{K^\omega(\tau^n \circ \phi \circ \sigma^n) \circ \eta^\omega_C} \circ \iota_1 \circ \alpha \\
  [\text{$\eta^\omega$ is natural}]
    &= \reducenext{\tilde\beta \circ \eta^\omega_C} \circ \tau^n \circ \phi \circ \sigma^n \circ \iota_1 \circ \alpha \\
  [\text{definition of $\tilde\beta$}]
    &= \alpha^{-1} \circ \reducenext{\beta \circ \tau^n \circ \phi \circ \sigma^n \circ \iota_1} \circ \alpha \\
  [\text{Lemma~\ref{kf.lem.calc-2}}]
    &= \alpha^{-1} \circ f^\alpha_{n+1} \circ \alpha = f_{n+1}.
  \end{align*}

  This shows that every $f_n$ belongs to the semigroup generated by $\tilde\alpha$, $\tilde\beta$, $\tilde\sigma$,
  $\tilde\tau$ and $\tilde\phi$, and we uniformly have that the length of the word representing $f_n$ is~$2n + 1$.
  Therefore, $\End(K^\omega(C))$ is strongly distorted and the \Sierpinski\ rank of $\End(K^\omega(C))$ is at most~5.
  Lemma~\ref{kf.lem.SSD} now yields that $\End(L)$ has the Bergman property if it is not finitely generated.
\end{PROOF}

\begin{COR}
    For the following \Fraisse\ limits $L$ we have that $\End(L)$ has the Bergman property:
    \begin{itemize}
    \item
      the random graph,
    \item
      the random digraph,
    \item
      the random poset,
    \item
      the rational Urysohn sphere (the \Fraisse\ limit of the category of all finite metric spaces with rational distances bounded by~1),
    \item
      the countable atomless Boolean algebra.
    \end{itemize}
\end{COR}

\begin{PROOF}
  It is easy to see that each of the categories involved has retractive natural (JEP). In the first four cases
  the existence of a retraction $r : K(L) \to L$ such that $r \circ \eta_L = \id_L$ where $L$ is the corresponding \Fraisse\ limit
  follows from the explicit construction of the \Katetov\ functor (Subsection~\ref{kf.subsec.ex}).

Let $(U,\varrho)$ denote the Urysohn sphere.
Note that each $p \in K(U)$ is {determined by} a finite set $F \subseteq U$ in the sense that
\begin{equation}\label{Eqernhguodg}
\varrho(p, u) = \min_{x \in F}\Bigl( \varrho(p,x) + \varrho(x,u) \Bigr)
\end{equation}
(see Section \ref{SecAppndxKattv} for more details, in particular, formula (\ref{EqKttvMppng})).
Note also that enlarging the set $F$, the equation above remains true, because of the triangle inequality.

Suppose that $U \subseteq X_0 \subseteq K(U)$ is such that $X_0 \setminus U$ is finite and a nonexpansive retraction $\map r{X_0}U$ has already been defined.
Fix $p \in K(U) \setminus X_0$.
Choose a finite set $F$ containing $r[X_0 \setminus U]$, such that (\ref{Eqernhguodg}) holds.
Let $A = (X_0 \setminus U) \cup F$.
Then $\map {r\restriction A}{A}{U}$ is a nonexpansive mapping (which is identity on $F$) therefore by \cite[Thm. 3.18]{KubInjRts} it extends to a nonexpansive mapping on $\map {\bar r}{A \cup \{p\}}U$.
Let $q = \bar r(p)$ and let $r' = \bar r \cup \id_U$.
We claim that $\map{r'}{X_0 \cup \{p\}}{U}$ is nonexpansive.

Fix $u \in U \setminus F$.
Using (\ref{Eqernhguodg}), we have $\varrho(p,u) = \varrho(p,s) + \varrho(s,u)$ for some $s \in F$.
Thus
$$\varrho(r'(p),r'(u)) = \varrho(\bar r(p), u) \le \varrho(\bar r(p), s) + \varrho(s, u) \le \varrho(p,s) + \varrho(s,u) = \varrho(p,u).$$
The last inequality follows from the fact that $\bar r$ is nonexpansive on $A$ and $\bar r(s) = s$.
This shows that $r'$ is a nonexpansive extension of $r$.
Easy induction shows the existence of a nonexpansive retraction of $K(U)$ onto $U$.

  Let us finally show that a retraction $\map r {K(L)}{L}$ also exists in case of the category of finite Boolean algebras.
  Recall from Example~\ref{kf.ex.BA} that for a finite Boolean algebra $B(A)$ whose set of atoms is $A$
  we have the \Katetov\ functor $K(B(A)) = B(\{0,1\} \times A)$ where $\eta_{B(A)} : B(A) \hookrightarrow B(\{0,1\} \times A)$ is
  the unique homomorphism which takes $a \in A$ to $\pair 0a \lor \pair 1a \in B(\{0,1\} \times A)$.
  It is now easy to see that for each finite Boolean algebra $B(A)$ there is a retraction $r_{B(A)} : K(B(A)) \to B(A)$
  which takes $\pair ia$ to $a$ $(i \in \{0, 1\})$
  and extends to the rest of $K(B(A))$ in an obvious way. Clearly, $r_{B(A)} \circ \eta_{B(A)} = \id_{B(A)}$.
  Let $L$ be the countable atomless Boolean algebra and let $B_1 \hookrightarrow B_2 \hookrightarrow \ldots$ be a chain of finite
  Boolean algebras whose colimit is $L$. Then the colimit of the chain
  $K(B_1) \hookrightarrow K(B_2) \hookrightarrow \ldots$ is $K(L)$ and the following diagram commutes:
  $$
    \XYMATRIX{
      & & & L \aremb@/_4mm/[ddd]_-{\eta_{L}}\\
      B_1 \aremb[r] \aremb@/_3mm/[d]_-{\eta_{B_1}} \aremb@/^3mm/[rrru] & B_2 \aremb[r] \aremb@/_3mm/[d]_-{\eta_{B_2}} \aremb@/^2mm/[rru] & \cdots \\
      K(B_1) \aremb[r] \ar@/_3mm/[u]_-{r_{B_1}} \aremb@/_3mm/[rrrd] & K(B_2) \aremb[r]  \ar@/_3mm/[u]_-{r_{B_1}} \aremb@/_2mm/[rrd] & \cdots\\
      & & & K(L) \ar@{.>}@/_4mm/[uuu]_-{r_{L}}
    }
  $$
  Since $K(L)$ is the colimit of the bottom chain, there is a unique mapping $r_L : K(L) \to L$ such that the diagram commutes.
  In particular, $r_L \circ \eta_L = \id_L$.
\end{PROOF}

\section{Appendix: the original \Katetov\ construction}\label{SecAppndxKattv}

For the sake of completeness we present the details of \Katetov's construction in the case of finite spaces.
Actually, we were unable to find any source were \Katetov's extensions of metric spaces are explicitly treated as a functor acting on nonexpansive mappings.

Given a metric space $X$ we shall denote its metric either by $\varrho$ or by $\varrho_X$.
Fix a finite metric space $X$ and denote by $K(X)$ the set of all functions $\map \phi X{[0,+\infty)}$ satisfying
$$|\phi(x_0) - \phi(x_1)| \le \varrho(x_0,x_1) \le \phi(x_0) + \phi(x_1)$$
for every $x_0, x_1 \in X$.
Elements of $K(X)$ are called \emph{\Katetov\ functions} on $X$.
Given $a \in X$, the function $\widehat a(x) = \varrho(x,a)$ is \Katetov, therefore it is natural to define $\map{\eta_X}X {K(X)}$ by $\eta_X(x) = \widehat x$.
Endow $K(X)$ with the metric
$$\varrho(\phi,\psi) = \max_{x \in X}|\phi(x) - \psi(x)|.$$
It is easy to see that $\eta_X$ is an isometric embedding.
Note that $K(X)$ is a Polish space, being a closed subspace of $\RR^X$.

We now fix a nonexpansive map $\map f X Y$ between finite metric spaces.
Given $\phi \in K(X)$, define
\begin{equation}\label{EqKttvMppng}
\phi^f(y) = \min_{x \in X}\Bigl(\varrho_Y(y,f(x)) + \phi(x)\Bigr).
\end{equation}

\begin{LEM}\label{Lmerbgibr}
$\phi^f \in K(Y)$ for every $\phi \in K(X)$. Furthermore, given $x \in X$, we have that $\phi^f(f(x)) \le \phi(x)$ and $\phi^f(f(x)) = \phi(x)$ whenever $f$ is an isometric embedding.
\end{LEM}

\begin{PROOF}
Fix $y_0, y_1 \in Y$ and assume $\phi^f(y_i) = \varrho_Y(y_i,f(x_i))+\phi(x_i)$ for $i=0,1$.
Then
\begin{align*}
\phi^f(y_0) &\le \varrho_Y(y_0,f(x_1)) + \phi(x_1) \\
&\le \varrho_Y(y_0,f(x_1)) - \varrho_Y(y_1,f(x_1)) + \varrho_Y(y_1,f(x_1)) + \phi(x_1) \\
&\le \varrho_Y(y_0,y_1) + \phi^f(y_1).
\end{align*}
Similarly, $\phi^f(y_1) \le \varrho_Y(y_0,y_1) + \phi^f(y_0)$.
Furthermore, using the fact that $f$ is nonexpansive and $\phi$ is \Katetov, we get
\begin{align*}
\varrho_Y(y_0,y_1) &\le \varrho_Y(y_0,f(x_0)) + \varrho_Y(f(x_0),f(x_1)) + \varrho_Y(y_1,f(x_1)) \\
&\le \varrho_Y(y_0,f(x_0)) + \varrho_X(x_0,x_1) + \varrho_Y(y_1,f(x_1)) \\
&\le \varrho_Y(y_0,f(x_0)) + \phi(x_0) + \phi(x_1) + \varrho_Y(y_1,f(x_1)) \\
&= \phi^f(y_0) + \phi^f(y_1).
\end{align*}
This shows that $\phi^f$ is a \Katetov\ function.
Inequality $\phi^f(f(x)) \le \phi(x)$ is trivial.
Finally, suppose $f$ is an isometric embedding and fix $x\in X$. Choose $x_1 \in X$ so that $\phi^f(f(x)) = \varrho_Y(f(x),f(x_1))+\phi(x_1)$.
Then
$$\phi^f(f(x)) = \varrho_X(x,x_1)+\phi(x_1) \ge \phi(x),$$
because $\phi$ is \Katetov.
This completes the proof.
\end{PROOF}

\begin{LEM}\label{Lmvbiuergsef}
$\varrho(\phi^f, \psi^f) \le \varrho(\phi, \psi)$ for every $\phi, \psi \in K(X)$.
Equality holds whenever $f$ is an isometric embedding.
\end{LEM}

\begin{PROOF}
Fix $y \in Y$.
Find $x_0 \in X$ such that $\phi^f(y) = \varrho(y,f(x_0)) + \phi(x_0)$.
Then $\psi^f(y) \le \varrho(y,f(x_0)) + \psi(x_0)$, therefore
$$\psi^f(y) - \phi^f(y) \le \psi(x_0) - \phi(x_0) \le |\psi(x_0)-\phi(x_0)| \le \varrho(\phi,\psi).$$
By symmetry, $\phi^f(y) - \psi^f(y) \le \varrho(\phi,\psi)$.
Thus $|\phi^f(y)-\psi^f(y)| \le \varrho(\phi,\psi)$ and hence $\varrho(\phi^f, \psi^f) \le \varrho(\phi, \psi)$.
Finally, if $f$ is an isometric embedding and $\varrho(\phi, \psi) = |\phi(x_0)-\psi(x_0)|$ then, using Lemma~\ref{Lmerbgibr}, we get
$$\varrho(\phi^f,\psi^f) \ge |\phi^f(x) - \psi^f(x)| = |\phi(x)-\psi(x)|$$
for every $x \in X$, which implies that $\varrho(\phi^f,\psi^f) \ge \varrho(\phi,\psi)$.
\end{PROOF}

\begin{LEM}
Given nonexpansive mappings $\map f X Y$, $\map g Y Z$ between finite metric spaces, it holds that $\phi^{g\cmp f} = (\phi^f)^g$ for every $\phi \in K(X)$.
\end{LEM}

\begin{PROOF}
Fix $z \in Z$.
We have
\begin{align*}
(\phi^f)^g(z) &= \min_{y \in Y} \Bigl( \varrho_Z(z,g(y))+\phi^f(y) \Bigr) \\
&= \min_{y \in Y,\; x\in X} \Bigl( \varrho_Z(z,g(y))+\varrho_Y(y,f(x))+\phi(x) \Bigr) \\
&\ge \min_{y \in Y,\; x\in X} \Bigl( \varrho_Z(z,g(y))+\varrho_Z(g(y),gf(x))+\phi(x) \Bigr) \\
&\ge \min_{x \in X} \Bigl( \varrho_Z(z,gf(x))+\phi(x) \Bigr) = \phi^{g\cmp f}(z).
\end{align*}
On the other hand, using Lemma~\ref{Lmerbgibr}, we get
\begin{align*}
(\phi^f)^g(z) &\le \min_{x\in X} \Bigl( \varrho(z,gf(x))+\phi^f(f(x)) \Bigr) \\
&\le \min_{x\in X} \Bigl( \varrho(z,gf(x))+\phi(x) \Bigr) = \phi^{g\cmp f}(z).
\end{align*}
\end{PROOF}

It is obvious that $\phi^{\id X} = \phi$, therefore defining
$$K(f)(\phi) = \phi^f$$
we obtain a covariant functor $K$ from the category of finite metric spaces into the category of Polish metric spaces, both considered with nonexpansive mappings.
Furthermore, $K$ preserves isometric embeddings (by the second part of Lemma~\ref{Lmvbiuergsef}).

\begin{LEM}
Given a nonexpansive mapping of finite metric spaces $\map f X Y$, the following diagram is commutative.
$$\xymatrix{
X \ar[d]_f \ar[rr]^{\eta_X} & & K(X) \ar[d]^{K(f)} \\
Y \ar[rr]_{\eta_Y} & & K(Y)
}$$
\end{LEM}

\begin{PROOF}
Fix $x \in X$.
We have $(K(f) \cmp \eta_X)(x) = K(f)(\widehat x) = (\widehat x)^f$ and
$(\eta_Y \cmp f)(x) = \eta_Y(f(x)) = \widehat{f(x)}$.
It remains to show that $(\widehat x)^f = \widehat{f(x)}$.
We have
\begin{align*}
(\widehat x)^f(y) &= \min_{t \in X}\Bigl(\varrho(y,f(t))+\varrho(x,t)\Bigr) \\
&\ge \min_{t \in X}\Bigl(\varrho(y,f(t))+\varrho(f(x),f(t))\Bigr) \ge \varrho(y,f(x)) = \widehat{f(x)}(y).
\end{align*}
On the other hand,
$$(\widehat x)^f(y) \le \varrho(y,f(x))+\varrho(x,x) = \varrho(y,f(x)) = \widehat{f(x)}(y).$$
Hence $(\widehat x)^f = \widehat{f(x)}$.
\end{PROOF}

The lemma above says that $\eta$ is a natural transformation from the identity functor into $K$.
The last fact just says that $K$ is a \Katetov\ functor.

\begin{PROP}
Let $\map e X Y$ be an isometric embedding such that $X$ is finite and $|Y \setminus X| = 1$.
Then there exists an isometric embedding $\map g Y {K(X)}$ such that $g \cmp e = \eta_X$.
\end{PROP}
\begin{PROOF}
We may assume that $Y = X \cup \{s\}$ and $e$ is the inclusion.
Let $\phi(x) = \varrho_Y(x,s)$.
Then $\phi$ is a \Katetov\ function on $X$ and hence, setting $g(s) = \phi$ and $g(x) = \widehat x$ for $x\in X$, we obtained the required embedding.
\end{PROOF}

Exactly the same arguments show the existence of a \Katetov\ functor for finite metric spaces with rational distances, leading to the rational Urysohn space.
{Finally, one can restrict the set of distances to the unit interval $[0,1]$ obtaining a \Katetov\ functor leading to the Urysohn sphere (or its rational variant).}
On the other hand, knowing that the category of finite metric spaces has one-point extension pushouts, Theorem~\ref{kf.thm.OEP} provides another \Katetov\ functor on the category of finite rational metric spaces.
The original \Katetov\ functor is better in the sense that, when working in the category of \emph{all} finite metric spaces, its values are complete separable metric spaces, which can be viewed as ``minimal" spaces realizing all one-point extensions.

\subsection{Conclusion}

As we have seen above, the original \Katetov\ construction deals with complete metric spaces, therefore it is formally out of the scope of our model-theoretic approach.
The same applies to the recent Ben Yaacov's construction~\cite{ben-yaacov} of a \Katetov\ functor on separable Banach spaces, leading to the so-called \emph{Gurari\u\i\ space}, the unique universal separable Banach space that is \emph{almost homogeneous}, namely, isometries between finite-dimensional subspaces can be approximated by bijective isometries of the entire space.
Both examples can be presented in the framework of continuous model theory \cite{ben-yaacov2}.
In the definition of a \Katetov\ functor one would need to relax the extension property, as the Gurari\u\i\ space satisfies only its approximate variant.

With some effort, one can adapt most of our arguments to categories of continuous models,
obtaining in particular the universality result of Uspenskij~\cite{uspenskij} as well as its counterpart concerning  monoids of nonexpansive mappings.
We have decided to present the theory of \Katetov\ functors in discrete model-theoretic setting in order to make it more clear and accessible.

It is possible to provide a purely category-theoretic framework for \Katetov\ functors.
Another direction is to study uncountable iterations of \Katetov\ functors, obtaining models of arbitrary cardinality that are homogeneous with respect to their finitely generated substructures.
This will be done elsewhere.

\section{Acknowledgements}

The authors would like to thank Igor Dolinka and Christian Pech for their helpful comments on the early version of this paper.


\begin{thebibliography}{99}

\bibitem{ben-yaacov2}
I.\ Ben Yaacov: \Fraisse\ limits of metric structures. J. Symb. Log. 80 (2015), no. 1, 100--115.

\bibitem{ben-yaacov}
  I.\ Ben Yaacov: The linear isometry group of the Gurarij space is universal.
 Proc. Amer. Math. Soc. 142 (2014), no. 7, 2459--2467

\bibitem{bilge-melleray}
  D.\ Bilge, J.\ Melleray: Elements of finite order in automorphism groups of homogeneous structures.
 Contrib. Discrete Math. 8 (2013), no. 2, 88--119.

\bibitem{bdd}
  A.\ Bonato, D.\ Deli\'c, I.\ Dolinka: All countable monoids embed into the monoid
  of the infinite random graph. Discrete Math.\ 310 (2010), 373--375.

\bibitem{CamNes}
P.\ Cameron, J.\ Ne\v set\v ril: Homomorphism-homogeneous relational structures. Combin. Probab. Comput. 15 (2006) 91--103.

\bibitem{dolinka}
  I.\ Dolinka: The Bergman property for endomorphism monoids of some \Fraisse\ limits.
  Forum Math. 26 (2014), no. 2, 357--376.

\bibitem{d-rp}
  I.\ Dolinka: The endomorphism monoid of the random poset contains all countable semigroups.
  Algebra Universalis\ 56 (2007), 469--474.

\bibitem{dolinka-masulovic}
  I.\ Dolinka, D.\ Ma\v sulovi\'c: A universality result for endomorphism monoids of some ultrahomogeneous structures.
  Proceedings of the Edinburgh Mathematical Society 55 (2012), 635--656.

\bibitem{gratzer-ua}
  G.\ Gr\"atzer: Universal algebra (2nd ed). Springer-Verlag, New York, 2008.


\bibitem{henson}
  C.\ W.\ Henson: A family of countable homogeneous graphs. Pac.\ J.\ Math.\ 38 (1971), 69--83.

\bibitem{herwig-lascar}
  B.\ Herwig, D.\ Lascar: Extending partial automorphisms and the profinite topology on free groups.
  Trans.\ Amer.\ Math.\ Soc.\ 352 (2000), 1985--2021.

\bibitem{katetov}
  M.\ Kat\v etov: On universal metric spaces. General topology and its relations
  to modern analysis and algebra. VI (Prague, 1986), Res.\ Exp.\ Math., vol.~16,
  Heldermann, Berlin, 1988, 323--330.

\bibitem{KePeTo}
A.\ S.\ Kechris, V.\ G.\ Pestov, S.\ Todorcevic:
\Fraisse\ limits, Ramsey theory, and topological dynamics of automorphism groups.
Geom. Funct. Anal. 15 (2005) 106--189.

\bibitem{KubInjRts} W.\ Kubi\'s: Injective objects and retracts of Fra\"\i ss\'e limits, Forum Math. 27 (2015) 807--842.

\bibitem{MacLane}
  S.\ Mac Lane: Categories for the working mathematician (2nd ed). Springer, 1978


\bibitem{maltcev-mitchell-ruskuc}
  V.\ Maltcev, J.\ D.\ Mitchell, N.\ Ru\v skuc: The Bergman property for semigroups. J.\ London Math.\ Soc.\ 80 (2009), 212--232.

\bibitem{peresse}
  Y.\ P\'eresse: Generating uncountable transformation semigroups. Ph.D.\ thesis, University of St Andrews, 2009.

\bibitem{solecki}
  S.\ Solecki: Notes on a strengthening of the Herwig-Lascar extension theorem. 2009. Unpublished
  note, available at \url{http://www.math.uiuc.edu/~ssolecki/papers/HervLascfin.pdf}

\bibitem{Urysohn}
{P.S.\ Urysohn}: {Sur un espace m\'etrique universel}, I, II,
Bull. Sci. Math. (2) 51 (1927), 43--64, 74--90.
\href{http://www.zentralblatt-math.org/zmath/en/advanced/?q=an:53.0556.01&format=complete}{JFM 53.0556.01}

\bibitem{uspenskij}
  V.\ V.\ Uspenskij: On the group of isometries of the Urysohn universal metric space.
  Commentat.\ Math.\ Univ.\ Carolinae 31 (1990), 181--182.

\end{thebibliography}
\end{document}